\documentclass[1pt]{amsart}

\usepackage{amsmath,amsthm, amscd, amssymb, amsfonts}

\usepackage{hhline}

\newcommand{\Ind}{\operatorname{Ind}}
\newcommand{\ord}{\operatorname{ord}}
\newcommand{\oper}{\#}

\DeclareMathOperator*{\prode}{\boxtimes}

\newcommand\toba{{\mathfrak B }}
\newcommand\trasp{\pi}
\newcommand{\gr}{\operatorname{gr}}
\newcommand{\trid}{\triangleright}

\newcommand\compvert[2]{\genfrac{}{}{-1pt}{0}{#1}{#2}}

\newcommand\mvert[2]{\begin{tiny}\begin{matrix}#1\vspace{-4pt}\\#2\end{matrix}\end{tiny}}
\DeclareMathOperator*{\Tim}{\times}
\newcommand{\Times}[2]{\sideset{_#1}{_#2}\Tim}

\newcommand{\cx}{{\daleth}}

\newcommand{\lu}{{\nu}}
\newcommand{\ld}{{\ell}}

\newcommand{\abofe}{{A(\lu, \ld)}}

\newcommand{\bbofe}{{B(\lu, \ld)}}

\newcommand{\proy}{\Pi}

\newcommand{\Lc}{{\mathcal L}}
\newcommand{\Y}{{\mathcal Y}}
\newcommand{\R}{{\mathcal R}}
\newcommand{\W}{{\mathcal W}}
\newcommand{\Zc}{{\mathcal Z}}
\newcommand{\daga}{{\dagger}}

\newcommand{\acts}{{\rightharpoondown}}

\newcommand{\ku}{\mathbb C}
\newcommand{\K}{{\mathcal K}}
\newcommand{\Z}{{\mathbb Z}}
\newcommand{\N}{{\mathbb N}}
\newcommand{\I}{{\mathbb I}}

\newcommand{\G}{{\mathbb G}}
\newcommand{\M}{{\mathcal M}}
\newcommand{\Q}{{\mathcal Q}}
\newcommand{\F}{{\mathcal F}}
\newcommand{\C}{{\mathcal C}}
\newcommand{\D}{{\mathcal D}}
\newcommand{\Ec}{{\mathbf E}}
\newcommand{\Ee}{{\mathcal E}}
\newcommand{\Kc}{{\mathbf K}}

\newcommand{\uno}{{\bf 1}}
\newcommand{\uv}{{\bf 1_v}}
\newcommand{\uh}{{\bf 1_h}}
\newcommand{\m}{\mathcal{M}}
\newcommand{\n}{\mathcal{N}}

\newcommand{\Dc}{{\mathbf D}}
\newcommand{\B}{{\mathcal B}}
\newcommand{\T}{{\mathcal T}}
\newcommand{\Hc}{{\mathcal H}}
\newcommand{\Vc}{{\mathcal V}}
\newcommand{\Pc}{{\mathcal P}}
\newcommand{\Oc}{{\mathcal O}}
\newcommand{\ydh}{{}^H_H\mathcal{YD}}

\newcommand{\Ss}{{\mathcal S}}
\newcommand{\Vect}{\operatorname{Vec}}
\newcommand{\End}{\operatorname{End}}
\newcommand{\Aut}{\operatorname{Aut}}
\newcommand{\Int}{\operatorname{Int}}
\newcommand{\Ext}{\operatorname{Ext}}
\newcommand\tr{\operatorname{tr}}
\newcommand\Rep{\operatorname{Rep}}

\newcommand\card{\operatorname{card}}
\newcommand{\unosigma}{{\bf 1}^{\sigma}}
\newcommand\tot{\operatorname{tot}}
\newcommand\sgn{\operatorname{sgn}}
\newcommand\ad{\operatorname{ad}}
\newcommand\Hom{\operatorname{Hom}}
\newcommand\opext{\operatorname{Opext}}
\newcommand\Tot{\operatorname{Tot}}
\newcommand\Map{\operatorname{Map}}
\newcommand{\fde}{{\triangleright}}
\newcommand{\fiz}{{\triangleleft}}
\newcommand{\wC}{\widehat{C}}
\newcommand{\la}{\langle}
\newcommand{\ra}{\rangle}
\newcommand{\Comod}{\mbox{\rm Comod\,}}
\newcommand{\Mod}{\mbox{\rm Mod\,}}
\newcommand{\Sg}{{\mathfrak S}}
\newcommand{\funcion}{\mathfrak p}
\newcommand{\tauo}{{\widehat\tau}}
\newcommand{\prin}{t}
\newcommand{\fin}{b}
\newcommand{\pri}{r}
\newcommand{\fine}{l}
\newcommand\rh{\sim_{H}}
\newcommand\rv{\sim_{V}}
\newcommand\rd{\sim_{D}}

\newcommand\V{\operatorname{Vec}}

\theoremstyle{plain}

\renewcommand{\baselinestretch}{1.2}

\title{The character tables of centralizers in  Weyl Group of $E_8$ IV}
\author{ \small Shouchuan Zhang,     Peng Wang,  Jing Cheng, Hui Yang
}
\address{ Mathematics Department of Hunan University,\newline \indent Changsha China,
410082, E-mail: z9491@yahoo.com.cn }
\date{}

\begin{document}
\newtheorem{Proposition}{\quad Proposition}[section]
\newtheorem{Theorem}{\quad Theorem}
\newtheorem{Definition}[Proposition]{\quad Definition}
\newtheorem{Corollary}[Proposition]{\quad Corollary}
\newtheorem{Lemma}[Proposition]{\quad Lemma}
\newtheorem{Example}[Proposition]{\quad Example}

\maketitle \addtocounter{section}{-1}

\numberwithin{equation}{section}

\date{}


\begin {abstract} To classify the finite dimensional pointed Hopf
algebras with Weyl group $G$ of $E_8$, we obtain  the
representatives of conjugacy classes of $G$ and  all character
tables of centralizers of these representatives by means of software
GAP. In this paper we only list character table 65--94.


keywords: GAP, Hopf algebra, Weyl group, character.
\end {abstract}

\section{Introduction}\label {s0}

This article is to contribute to the classification of
finite-dimensional complex pointed Hopf algebras  with Weyl groups
of $E_8$.
 Many papers are about the classification of finite dimensional
pointed Hopf algebras, for example,  \cite{AS98b, AS02, AS00, AS05,
He06, AHS08, AG03, AFZ08, AZ07,Gr00,  Fa07, AF06, AF07, ZZC04, ZCZ08}.

In these research  ones need  the centralizers and character tables
of groups. In this paper we obtain   the representatives of
conjugacy classes of Weyl groups of $E_8$  and all
character tables of centralizers of these representatives by means
of software GAP. In this paper we only list character table  65--94.

By the Cartan-Killing classification of simple Lie algebras over complex field the Weyl groups
to
be considered are $W(A_l), (l\geq 1); $
 $W(B_l), (l \geq 2);$
 $W(C_l), (l \geq 2); $
 $W(D_l), (l\geq 4); $
 $W(E_l), (8 \geq l \geq 6); $
 $W(F_l), ( l=4 );$
$W(G_l), (l=2)$.

It is otherwise desirable to do this in view of the importance of Weyl groups in the theories of
Lie groups, Lie algebras and algebraic groups. For example, the irreducible representations of
Weyl groups were obtained by Frobenius, Schur and Young. The conjugace classes of  $W(F_4)$
were obtained by Wall
 \cite {Wa63} and its character tables were obtained by Kondo \cite {Ko65}.
The conjugace classes  and  character tables of $W(E_6),$ $W(E_7)$ and $W(E_8)$ were obtained
by  Frame \cite {Fr51}.
Carter gave a unified description of the conjugace classes of
 Weyl groups of simple Lie algebras  \cite  {Ca72}.

\section {Program}

 By using the following
program in GAP,  we  obtain  the representatives of conjugacy
classes of Weyl groups of $E_6$  and all character tables of
centralizers of these representatives.

gap$>$ L:=SimpleLieAlgebra("E",6,Rationals);;

gap$>$ R:=RootSystem(L);;

 gap$>$ W:=WeylGroup(R);Display(Order(W));

gap $>$ ccl:=ConjugacyClasses(W);;

gap$>$ q:=NrConjugacyClasses(W);; Display (q);

gap$>$ for i in [1..q] do

$>$ r:=Order(Representative(ccl[i]));Display(r);;

$>$ od; gap

$>$ s1:=Representative(ccl[1]);;cen1:=Centralizer(W,s1);;

gap$>$ cl1:=ConjugacyClasses(cen1);

gap$>$ s1:=Representative(ccl[2]);;cen1:=Centralizer(W,s1);;

 gap$>$
cl2:=ConjugacyClasses(cen1);

gap$>$ s1:=Representative(ccl[3]);;cen1:=Centralizer(W,s1);;

gap$>$ cl3:=ConjugacyClasses(cen1);

 gap$>$
s1:=Representative(ccl[4]);;cen1:=Centralizer(W,s1);;

 gap$>$
cl4:=ConjugacyClasses(cen1);

 gap$>$
s1:=Representative(ccl[5]);;cen1:=Centralizer(W,s1);;

gap$>$ cl5:=ConjugacyClasses(cen1);

 gap$>$
s1:=Representative(ccl[6]);;cen1:=Centralizer(W,s1);;

gap$>$ cl6:=ConjugacyClasses(cen1);

gap$>$ s1:=Representative(ccl[7]);;cen1:=Centralizer(W,s1);;

gap$>$ cl7:=ConjugacyClasses(cen1);

gap$>$ s1:=Representative(ccl[8]);;cen1:=Centralizer(W,s1);;

gap$>$ cl8:=ConjugacyClasses(cen1);

gap$>$ s1:=Representative(ccl[9]);;cen1:=Centralizer(W,s1);;

gap$>$ cl9:=ConjugacyClasses(cen1);

gap$>$ s1:=Representative(ccl[10]);;cen1:=Centralizer(W,s1);;

 gap$>$
cl10:=ConjugacyClasses(cen1);

gap$>$ s1:=Representative(ccl[11]);;cen1:=Centralizer(W,s1);;

gap$>$ cl11:=ConjugacyClasses(cen1);

gap$>$ s1:=Representative(ccl[12]);;cen1:=Centralizer(W,s1);;

gap$>$ cl2:=ConjugacyClasses(cen1);

gap$>$ s1:=Representative(ccl[13]);;cen1:=Centralizer(W,s1);;

 gap$>$
cl13:=ConjugacyClasses(cen1);

 gap$>$
s1:=Representative(ccl[14]);;cen1:=Centralizer(W,s1);;

gap$>$ cl14:=ConjugacyClasses(cen1);

gap$>$ s1:=Representative(ccl[15]);;cen1:=Centralizer(W,s1);;

gap$>$ cl15:=ConjugacyClasses(cen1);

gap$>$ s1:=Representative(ccl[16]);;cen1:=Centralizer(W,s1);;

gap$>$ cl16:=ConjugacyClasses(cen1);

gap$>$ s1:=Representative(ccl[17]);;cen1:=Centralizer(W,s1);;

gap$>$ cl17:=ConjugacyClasses(cen1);

gap$>$ s1:=Representative(ccl[18]);;cen1:=Centralizer(W,s1);;

gap$>$ cl18:=ConjugacyClasses(cen1);

gap$>$ s1:=Representative(ccl[19]);;cen1:=Centralizer(W,s1);;

gap$>$ cl19:=ConjugacyClasses(cen1);

gap$>$ s1:=Representative(ccl[20]);;cen1:=Centralizer(W,s1);;

gap$>$ cl20:=ConjugacyClasses(cen1);

gap$>$ s1:=Representative(ccl[21]);;cen1:=Centralizer(W,s1);;

 gap$>$
cl21:=ConjugacyClasses(cen1);

gap$>$ s1:=Representative(ccl[22]);;cen1:=Centralizer(W,s1);;

 gap$>$
cl22:=ConjugacyClasses(cen1);

gap$>$ s1:=Representative(ccl[23]);;cen1:=Centralizer(W,s1);;

gap$>$ cl23:=ConjugacyClasses(cen1);

gap$>$ s1:=Representative(ccl[24]);;cen1:=Centralizer(W,s1);;

$>$ cl24:=ConjugacyClasses(cen1);

gap$>$ s1:=Representative(ccl[25]);;cen1:=Centralizer(W,s1);

gap$>$ cl25:=ConjugacyClasses(cen1);

gap$>$ for i in [1..q] do

$>$ s:=Representative(ccl[i]);;cen:=Centralizer(W,s);;

$>$ char:=CharacterTable(cen);;Display (cen);Display(char);

 $>$ od;

The programs for Weyl groups of $E_7$, $E_8$, $F_4$ and $ G_2$ are
similar.

{\small

\section {$E_8$}

The generators of $G^{s_{65}}$ are:\\
$\left(
\begin{array}{cccccccc}0&-2&0&1&0&0&0&-1\\1&-2&-1&1&1&-1&0&
-1\\1&-3&0&1&1&-1&0&
-2\\ 2&-4&-1&1&2&-1&0&-3\\ 2&-3&-1&1&1&0&-1&-2 \\
      2&-3&-1&1&1&0&-1&-1\\1&-2&-1&1&1&0&-1&-1 \\
      1&-1&-1&1&0&0&0&-1\end{array}\right)$
  $\left(\begin{array}{cccccccc}0&-1&0&1&0&0&-1&1\\ -1&0&0&1&0&0&-1&1 \\ -1&-1&1&1&0&0&-1&1
      \\ -1&-1&0&2&0&0&-1&1 \\ -1&-1&0&1&1&0&-1&1\\-1&-1&0&1&0&1&-1&1
      \\ -1&-1&0&1&0&0&0&1\\ 0&0&0&0&0&0&0&1\end{array}\right)$
  $\left(\begin{array}{cccccccc} 1&-1&0&0&1&-1&0&1\\0&-1&0&0&2&-2&0&2\\0&-2&1&0&2&-2&0&2
     \\ 0&-3&0&1&3&-3&0&3\\0&-2&0&0&3&-2&0&2\\0&-2&0&0&2&-1&0&2
     \\ 0&-1&0&0&1&-1&1&1\\0&0&0&0&0&0&0&1 \end{array}\right)$.

The representatives of conjugacy classes of   $G^{s_{65}}$ are:\\
$\left(\begin{array}{cccccccc} 1&0&0&0&0&0&0&0 \\ 0&1&0&0&
0&0&0&0 \\ 0&0&1&0&0&0&0&0 \\
        0&0&0&1&0&0&0&0 \\ 0&0&0&0&1&0&0&0 \\ 0&0&0&0&0&1&0&0 \\
        0&0&0&0&0&0&1&0 \\ 0&0&0&0&0&0&0&1 \end{array}\right)$,
  $\left(\begin{array}{cccccccc} -2&-1&1&1&-1&0&1&-1 \\ -2&-1&0&2&-2&1&1&-2 \\
        -3&-1&1&2&-2&0&2&-3 \\ -4&-2&1&3&-3&1&3&-5 \\
        -3&-2&1&2&-2&1&2&-4 \\ -2&-2&1&2&-2&1&1&-3 \\
        -1&-2&0&2&-1&0&1&-2 \\ 0&-1&0&1&-1&0&1&-1 \end{array}\right)$,
  $\left(\begin{array}{cccccccc} -2&0&0&0&0&1&0&0 \\ -3&-1&1&0&0&1&0&0 \\ -3&-1&0&0&0&2&0&0 \\
        -5&-1&1&-1&1&2&0&0 \\ -4&0&1&-1&0&2&0&0 \\ -3&0&0&0&0&1&0&0 \\
        -2&0&0&0&0&1&0&-1 \\ -1&0&0&0&0&0&1&-1 \end{array}\right)$,
  $\left(\begin{array}{cccccccc} -2&0&0&0&0&1&0&0 \\ -3&0&1&0&-1&2&0&-1 \\ -3&0&0&0&-1&3&0&-1 \\
        -5&1&1&-1&-1&4&0&-2 \\ -4&1&1&-1&-1&3&0&-1 \\
        -3&1&0&0&-1&2&0&-1 \\ -2&0&0&0&0&1&0&-1 \\ -1&0&0&0&0&0&1&-1 \end{array}\right)$,
  $\left(\begin{array}{cccccccc} -2&0&1&0&0&-1&1&1 \\ -3&-1&1&1&-1&0&1&1 \\ -4&0&1&1&-1&-1&2&1 \\
        -6&-1&2&2&-2&-1&3&1 \\ -5&-1&2&1&-1&-1&3&0 \\
        -4&-1&2&1&-1&-1&2&0 \\ -2&-1&1&1&-1&-1&2&0 \\
        -1&0&1&0&-1&0&1&0 \end{array}\right)$,
  $\left(\begin{array}{cccccccc} -2&0&1&0&0&-1&1&1 \\ -2&-1&1&0&0&-1&2&1 \\ -3&0&1&0&0&-2&3&1 \\
        -4&-1&2&0&0&-3&5&1 \\ -4&-1&2&0&0&-2&4&0 \\ -3&-1&2&0&0&-2&3&0 \\
        -2&-1&1&1&-1&-1&2&0 \\ -1&0&1&0&-1&0&1&0 \end{array}\right)$,
  $\left(\begin{array}{cccccccc} -2&0&2&-1&0&1&-1&1 \\ -3&0&2&-1&0&1&-1&2 \\ -4&0&4&-2&0&1&-1&2 \\
        -5&0&5&-3&0&2&-2&4 \\ -4&0&4&-2&-1&2&-1&3 \\
        -3&-1&3&-1&-1&2&-1&2 \\ -2&0&2&-1&-1&2&-1&1 \\
        -1&0&1&-1&0&1&0&0 \end{array}\right)$,
  $\left(\begin{array}{cccccccc} -2&0&2&-1&0&1&-1&1 \\ -3&1&2&-1&-1&2&-1&1 \\ -4&1&4&-2&-1&2&-1&1 \\
        -5&2&5&-3&-2&4&-2&2 \\ -4&1&4&-2&-2&3&-1&2 \\
        -3&0&3&-1&-2&3&-1&1 \\ -2&0&2&-1&-1&2&-1&1 \\
        -1&0&1&-1&0&1&0&0 \end{array}\right)$,
  $\left(\begin{array}{cccccccc} -2&1&1&0&-1&0&1&0 \\ -3&0&1&1&-2&1&1&0 \\ -4&1&1&1&-2&0&2&0 \\
        -6&0&2&2&-3&0&3&0 \\ -5&0&2&1&-2&0&3&-1 \\ -4&0&2&1&-2&0&2&-1 \\
        -2&0&1&1&-2&0&2&-1 \\ -1&0&1&0&-1&0&1&0 \end{array}\right)$,
  $\left(\begin{array}{cccccccc} -1&-1&0&0&1&0&0&1 \\ -3&-1&1&0&1&0&0&1 \\ -3&-2&1&0&1&1&0&1 \\
        -5&-2&1&0&2&1&0&1 \\ -4&-1&1&-1&2&1&0&1 \\ -3&-1&0&0&1&1&0&1 \\
        -2&-1&0&0&1&0&1&0 \\ -1&0&0&0&0&0&1&0 \end{array}\right)$,
  $\left(\begin{array}{cccccccc} -1&-1&0&0&1&0&0&1 \\ -2&-1&1&-1&2&-1&1&1 \\ -2&-2&1&-1&2&0&1&1 \\
        -3&-2&1&-2&4&-1&2&1 \\ -3&-1&1&-2&3&0&1&1 \\
        -2&-1&0&-1&2&0&1&1 \\ -2&-1&0&0&1&0&1&0 \\ -1&0&0&0&0&0&1&0 \end{array}\right)$,
 $\left(\begin{array}{cccccccc} -1&-1&1&0&-1&1&0&1 \\ -2&-1&1&1&-2&1&0&1 \\
        -2&-2&2&1&-3&2&0&1 \\ -3&-2&2&2&-5&3&0&2 \\ -2&-2&1&2&-4&2&1&1 \\
        -1&-2&0&2&-3&2&0&1 \\ 0&-1&0&1&-2&1&0&1 \\ 0&0&0&0&-1&1&0&0 \end{array}\right)$,
 $\left(\begin{array}{cccccccc} -1&-1&1&0&-1&1&0&1 \\ -1&-1&1&0&-1&0&1&1 \\
        -1&-2&2&0&-2&1&1&1 \\ -1&-2&2&0&-3&1&2&2 \\ -1&-2&1&1&-3&1&2&1 \\
        0&-2&0&1&-2&1&1&1 \\ 0&-1&0&1&-2&1&0&1 \\ 0&0&0&0&-1&1&0&0 \end{array}\right)$,
  $\left(\begin{array}{cccccccc} -1&0&0&0&0&0&0&0 \\ -1&-1&0&1&-1&1&-1&0 \\ -1&0&-1&1&-1&1&-1&0 \\
        -2&0&0&1&-2&2&-2&0 \\ -1&0&0&1&-2&1&-1&0 \\ -1&0&0&1&-1&0&-1&0 \\
        0&0&0&0&0&0&-1&0 \\ 0&0&0&0&0&0&0&-1 \end{array}\right)$,
  $\left(\begin{array}{cccccccc} -1&0&0&0&0&0&0&0 \\ 0&-1&0&0&0&0&0&0 \\ 0&0&-1&0&0&0&0&0 \\
        0&0&0&-1&0&0&0&0 \\ 0&0&0&0&-1&0&0&0 \\ 0&0&0&0&0&-1&0&0 \\
        0&0&0&0&0&0&-1&0 \\ 0&0&0&0&0&0&0&-1 \end{array}\right)$,
  $\left(\begin{array}{cccccccc} -1&0&0&0&0&1&0&0 \\ -3&0&1&0&0&1&0&0 \\ -3&-1&1&0&0&2&0&0 \\
        -5&-1&1&0&1&2&0&0 \\ -4&0&1&-1&1&2&0&0 \\ -3&0&0&0&0&2&0&0 \\
        -2&0&0&0&0&1&1&-1 \\ -1&0&0&0&0&0&1&0 \end{array}\right)$,
  $\left(\begin{array}{cccccccc} -1&0&1&-1&2&-1&-1&1 \\ -2&0&1&-1&2&0&-2&2 \\ -3&1&2&-2&3&-1&-2&2 \\
        -4&1&3&-3&4&0&-4&3 \\ -3&1&3&-3&3&0&-3&3 \\ -3&0&3&-2&2&0&-2&2 \\
        -2&0&2&-2&2&0&-1&1 \\ -1&0&1&-1&1&0&0&0 \end{array}\right)$,
  $\left(\begin{array}{cccccccc} -1&0&1&-1&2&-1&-1&1 \\ -1&0&1&-2&3&-1&-1&2 \\
        -2&1&2&-3&4&-2&-1&2 \\ -2&1&3&-5&6&-2&-2&3 \\
        -2&1&3&-4&4&-1&-2&3 \\ -2&0&3&-3&3&-1&-1&2 \\
        -2&0&2&-2&2&0&-1&1 \\ -1&0&1&-1&1&0&0&0 \end{array}\right)$,
  $\left(\begin{array}{cccccccc} -1&0&1&0&-2&2&0&0 \\ -2&0&1&1&-3&2&0&0 \\ -2&-1&2&1&-4&3&0&0 \\
        -3&-1&2&2&-6&4&0&1 \\ -2&-1&1&2&-5&3&1&0 \\ -1&-1&0&2&-4&3&0&0 \\
        0&0&0&1&-3&2&0&0 \\ 0&0&0&0&-1&1&0&0 \end{array}\right)$,
  $\left(\begin{array}{cccccccc} -1&0&1&0&-1&0&2&-2 \\ -1&0&0&1&-2&1&2&-3 \\ -2&0&1&1&-2&0&3&-4 \\
        -3&-1&1&2&-3&1&4&-6 \\ -2&-1&1&1&-2&1&3&-5 \\
        -1&-1&1&1&-2&1&2&-4 \\ 0&-1&0&1&-1&0&2&-3 \\
        0&-1&0&1&-1&0&1&-1 \end{array}\right)$,
  $\left(\begin{array}{cccccccc} -1&1&0&-1&0&1&1&-1 \\ -2&1&1&-1&-1&2&1&-2 \\ -2&1&0&-1&-1&3&1&-2 \\
        -4&2&1&-2&-1&4&1&-3 \\ -3&2&1&-2&-1&3&1&-2 \\
        -2&2&0&-1&-1&2&1&-2 \\ -1&1&0&-1&0&1&1&-2 \\ -1&0&0&0&0&0&1&-1 \end{array}\right)$,
  $\left(\begin{array}{cccccccc} -1&1&0&0&-1&1&0&-1 \\ -1&0&0&1&-2&2&-1&-1 \\
        -1&1&-1&1&-2&2&-1&-1 \\ -2&1&0&1&-3&3&-2&-1 \\
        -1&1&0&1&-3&2&-1&-1 \\ -1&1&0&1&-2&1&-1&-1 \\
        0&1&0&0&-1&1&-1&-1 \\ 0&0&0&0&0&0&0&-1 \end{array}\right)$,
  $\left(\begin{array}{cccccccc} -1&1&0&0&0&0&-1&2 \\ -2&0&1&0&0&0&-1&3 \\ -2&1&0&0&0&0&-1&4 \\
        -3&1&1&0&0&-1&-1&6 \\ -3&1&1&0&0&-1&0&4 \\ -3&1&1&0&0&-1&0&3 \\
        -2&1&1&0&-1&0&0&2 \\ -1&1&1&-1&0&0&0&1 \end{array}\right)$,
  $\left(\begin{array}{cccccccc} -1&1&0&0&0&0&-1&2 \\ -2&1&1&0&-1&1&-1&2 \\ -2&2&0&0&-1&1&-1&3 \\
        -3&3&1&0&-2&1&-1&4 \\ -3&2&1&0&-1&0&0&3 \\ -3&2&1&0&-1&0&0&2 \\
        -2&1&1&0&-1&0&0&2 \\ -1&1&1&-1&0&0&0&1 \end{array}\right)$,
  $\left(\begin{array}{cccccccc} -1&1&1&-1&1&0&-1&0 \\ -2&1&1&-1&1&1&-2&1 \\ -3&2&2&-2&2&0&-2&1 \\
        -4&2&3&-3&3&1&-4&2 \\ -3&2&3&-3&2&1&-3&2 \\ -3&1&3&-2&1&1&-2&1 \\
        -2&1&2&-2&1&1&-1&0 \\ -1&0&1&-1&1&0&0&0 \end{array}\right)$,
  $\left(\begin{array}{cccccccc} -1&1&2&-2&0&1&0&0 \\ -2&2&2&-2&-1&2&0&0 \\ -3&2&4&-3&-1&2&0&0 \\
        -4&3&5&-4&-2&4&-1&1 \\ -3&2&4&-3&-2&3&0&1 \\
        -2&1&3&-2&-2&3&0&0 \\ -1&1&2&-2&-1&2&0&0 \\ -1&0&1&-1&0&1&0&0 \end{array}\right)$,
 $\left(\begin{array}{cccccccc} 0&-2&0&1&0&0&0&-1 \\ 0&-2&-1&2&0&0&-1&-1 \\
        0&-3&0&2&0&0&-1&-2 \\ 0&-4&-1&3&0&1&-2&-3 \\
        1&-3&-1&2&0&1&-2&-2 \\ 1&-3&-1&2&0&1&-2&-1 \\
        1&-2&-1&1&1&0&-1&-1 \\ 1&-1&-1&1&0&0&0&-1 \end{array}\right)$,
  $\left(\begin{array}{cccccccc} 0&-2&0&1&0&0&0&-1 \\ 1&-2&-1&1&1&-1&0&-1 \\ 1&-3&0&1&1&-1&0&-2 \\
        2&-4&-1&1&2&-1&0&-3 \\ 2&-3&-1&1&1&0&-1&-2 \\
        2&-3&-1&1&1&0&-1&-1 \\ 1&-2&-1&1&1&0&-1&-1 \\
        1&-1&-1&1&0&0&0&-1 \end{array}\right)$,
  $\left(\begin{array}{cccccccc} 0&-1&-1&2&-2&1&0&0 \\ 0&-2&-1&3&-2&0&0&0 \\ 1&-3&-2&4&-3&1&0&0 \\
        1&-4&-3&6&-4&0&1&0 \\ 1&-3&-3&5&-3&0&1&-1 \\
        1&-2&-3&4&-2&0&0&0 \\ 1&-1&-2&3&-2&0&0&0 \\ 1&0&-1&1&-1&0&0&0 \end{array}\right)$,
  $\left(\begin{array}{cccccccc} 0&-1&-1&2&-2&1&0&0 \\ 0&-1&-1&3&-3&1&0&-1 \\
        1&-2&-2&4&-4&2&0&-1 \\ 1&-2&-3&6&-6&2&1&-2 \\
        1&-2&-3&5&-4&1&1&-2 \\ 1&-1&-3&4&-3&1&0&-1 \\
        1&-1&-2&3&-2&0&0&0 \\ 1&0&-1&1&-1&0&0&0 \end{array}\right)$,
  $\left(\begin{array}{cccccccc} 0&-1&0&1&-1&1&0&-2 \\ 0&-1&-1&2&-1&1&-1&-2 \\
        0&-2&0&2&-1&1&-1&-3 \\ 0&-3&-1&3&-1&2&-2&-4 \\
        1&-2&-1&2&-1&2&-2&-3 \\ 1&-2&-1&2&-1&2&-2&-2 \\
        1&-1&-1&1&0&1&-1&-2 \\ 1&-1&-1&1&0&0&0&-1 \end{array}\right)$,
  $\left(\begin{array}{cccccccc} 0&-1&0&1&0&0&-1&1 \\ -1&-1&0&1&1&-1&-1&2 \\ -1&-2&1&1&1&-1&-1&2 \\
        -1&-3&0&2&2&-2&-1&3 \\ -1&-2&0&1&2&-1&-1&2 \\
        -1&-2&0&1&1&0&-1&2 \\ -1&-1&0&1&0&0&0&1 \\ 0&0&0&0&0&0&0&1 \end{array}\right)$,
  $\left(\begin{array}{cccccccc} 0&-1&0&1&0&0&-1&1 \\ -1&0&0&1&0&0&-1&1 \\ -1&-1&1&1&0&0&-1&1 \\
        -1&-1&0&2&0&0&-1&1 \\ -1&-1&0&1&1&0&-1&1 \\ -1&-1&0&1&0&1&-1&1 \\
        -1&-1&0&1&0&0&0&1 \\ 0&0&0&0&0&0&0&1 \end{array}\right)$,
  $\left(\begin{array}{cccccccc} 0&0&-1&1&1&-1&-1&0 \\ 0&-1&-1&1&2&-1&-2&1 \\ 0&0&-2&1&3&-2&-2&1 \\
        0&-1&-2&1&5&-3&-3&1 \\ 0&0&-1&0&4&-2&-3&1 \\
        -1&0&0&0&3&-2&-2&1 \\ -1&0&0&0&2&-1&-1&0 \\ 0&0&0&0&1&-1&0&0 \end{array}\right)$,
  $\left(\begin{array}{cccccccc} 0&0&-1&1&1&-1&-1&0 \\ 0&0&-1&1&1&0&-2&0 \\ 0&1&-2&1&2&-1&-2&0 \\
        0&1&-2&1&3&-1&-3&-1 \\ 0&1&-1&0&3&-1&-3&0 \\
        -1&1&0&0&2&-1&-2&0 \\ -1&0&0&0&2&-1&-1&0 \\ 0&0&0&0&1&-1&0&0 \end{array}\right)$,
  $\left(\begin{array}{cccccccc} 0&0&0&1&-1&0&-1&0 \\ 1&-1&-1&2&-1&0&-2&1 \\
        1&0&-1&2&-1&-1&-2&1 \\ 2&-1&-1&3&-2&-1&-3&2 \\
        2&-1&-1&3&-2&-1&-2&1 \\ 1&-1&0&2&-1&-1&-2&1 \\
        1&0&0&1&-1&0&-2&1 \\ 1&0&0&0&0&0&-1&0 \end{array}\right)$,
  $\left(\begin{array}{cccccccc} 0&0&0&1&-1&0&-1&0 \\ 1&0&-1&2&-2&1&-2&0 \\ 1&1&-1&2&-2&0&-2&0 \\
        2&1&-1&3&-4&1&-3&0 \\ 2&0&-1&3&-3&0&-2&0 \\ 1&0&0&2&-2&0&-2&0 \\
        1&0&0&1&-1&0&-2&1 \\ 1&0&0&0&0&0&-1&0 \end{array}\right)$,
  $\left(\begin{array}{cccccccc} 0&2&0&-1&0&0&0&1 \\ -1&2&1&-1&-1&1&0&1 \\ -1&3&0&-1&-1&1&0&2 \\
        -2&4&1&-1&-2&1&0&3 \\ -2&3&1&-1&-1&0&1&2 \\ -2&3&1&-1&-1&0&1&1 \\
        -1&2&1&-1&-1&0&1&1 \\ -1&1&1&-1&0&0&0&1 \end{array}\right)$,
  $\left(\begin{array}{cccccccc} 1&-1&-2&2&0&-1&0&0 \\ 1&-2&-2&3&0&-1&-1&0 \\ 2&-2&-4&4&0&-1&-1&0 \\
        2&-3&-5&6&0&-2&-1&-1 \\ 2&-2&-4&4&1&-2&-1&-1 \\
        1&-1&-3&3&1&-2&-1&0 \\ 1&-1&-2&2&1&-2&0&0 \\ 1&0&-1&1&0&-1&0&0 \end{array}\right)$,
  $\left(\begin{array}{cccccccc} 1&-1&-2&2&0&-1&0&0 \\ 2&-2&-2&2&1&-2&0&0 \\
        3&-2&-4&3&1&-2&0&0 \\ 4&-3&-5&4&2&-4&1&-1 \\
        3&-2&-4&3&2&-3&0&-1 \\ 2&-1&-3&2&2&-3&0&0 \\
        1&-1&-2&2&1&-2&0&0 \\ 1&0&-1&1&0&-1&0&0 \end{array}\right)$,
  $\left(\begin{array}{cccccccc} 1&-1&-1&1&0&1&-2&0 \\ 1&-1&-1&1&1&0&-3&1 \\ 2&-2&-1&1&1&1&-4&1 \\
        3&-2&-2&1&2&1&-6&2 \\ 3&-1&-2&1&1&1&-5&2 \\ 2&-1&-2&1&1&1&-4&2 \\
        1&0&-1&0&1&1&-3&1 \\ 1&0&-1&0&1&0&-1&0 \end{array}\right)$,
  $\left(\begin{array}{cccccccc} 1&-1&-1&1&0&1&-2&0 \\ 1&0&-1&1&0&1&-3&0 \\ 2&-1&-1&1&0&2&-4&0 \\
        3&0&-2&1&0&3&-6&0 \\ 3&0&-2&1&0&2&-5&1 \\ 2&0&-2&1&0&2&-4&1 \\
        1&0&-1&0&1&1&-3&1 \\ 1&0&-1&0&1&0&-1&0 \end{array}\right)$,
  $\left(\begin{array}{cccccccc} 1&-1&0&1&0&-1&-1&1 \\ 1&-1&-1&2&0&-1&-2&2 \\ 1&-1&0&2&0&-2&-2&2 \\
        2&-2&-1&4&-1&-2&-3&3 \\ 2&-2&-1&3&0&-2&-2&2 \\
        1&-2&0&2&0&-1&-2&2 \\ 1&-1&0&1&0&-1&-1&2 \\ 1&0&0&0&0&0&-1&1 \end{array}\right)$,
  $\left(\begin{array}{cccccccc} 1&-1&0&1&0&-1&-1&1 \\ 2&-1&-1&1&1&-2&-1&2 \\
        2&-1&0&1&1&-3&-1&2 \\ 4&-2&-1&2&1&-4&-1&3 \\
        3&-2&-1&2&1&-3&-1&2 \\ 2&-2&0&1&1&-2&-1&2 \\
        1&-1&0&1&0&-1&-1&2 \\ 1&0&0&0&0&0&-1&1 \end{array}\right)$,
  $\left(\begin{array}{cccccccc} 1&0&-2&2&-1&0&0&-1 \\ 1&-1&-2&3&-1&0&-1&-1 \\
        2&-1&-4&4&-1&0&-1&-1 \\ 2&-2&-5&6&-1&-1&-1&-2 \\
        2&-1&-4&4&0&-1&-1&-2 \\ 1&0&-3&3&0&-1&-1&-1 \\
        1&0&-2&2&0&-1&0&-1 \\ 1&0&-1&1&0&-1&0&0 \end{array}\right)$,
  $\left(\begin{array}{cccccccc} 1&0&-1&0&1&0&-2&2 \\ 0&0&0&0&1&0&-3&3 \\ 1&0&-1&0&1&1&-4&4 \\
        1&1&-1&0&1&1&-6&6 \\ 1&1&-1&0&1&0&-4&5 \\ 0&1&-1&0&1&0&-3&4 \\
        0&1&0&-1&1&0&-2&3 \\ 0&1&0&-1&1&0&-1&1 \end{array}\right)$,
  $\left(\begin{array}{cccccccc} 1&0&-1&0&1&0&-2&2 \\ 1&0&0&-1&2&-1&-2&3 \\ 2&0&-1&-1&2&0&-3&4 \\
        3&1&-1&-2&3&-1&-4&6 \\ 2&1&-1&-1&2&-1&-3&5 \\
        1&1&-1&-1&2&-1&-2&4 \\ 0&1&0&-1&1&0&-2&3 \\ 0&1&0&-1&1&0&-1&1 \end{array}\right)$,
  $\left(\begin{array}{cccccccc} 1&0&-1&1&-2&1&1&-1 \\ 1&0&-1&2&-3&1&1&-2 \\
        2&-1&-2&3&-4&2&1&-2 \\ 2&-1&-3&5&-6&2&2&-3 \\
        2&-1&-3&4&-4&1&2&-3 \\ 2&0&-3&3&-3&1&1&-2 \\
        2&0&-2&2&-2&0&1&-1 \\ 1&0&-1&1&-1&0&0&0 \end{array}\right)$,
  $\left(\begin{array}{cccccccc} -2&-1&1&1&-1&0&1&-1 \\ -2&-2&0&2&-1&0&1&-1 \\
        -3&-2&1&2&-1&-1&2&-2 \\ -4&-4&1&3&-1&-1&3&-3 \\
        -3&-3&1&2&-1&0&2&-3 \\ -2&-3&1&2&-1&0&1&-2 \\
        -1&-2&0&2&-1&0&1&-2 \\ 0&-1&0&1&-1&0&1&-1 \end{array}\right)$,
  $\left(\begin{array}{cccccccc} 1&0&0&1&-1&0&-1&0 \\ 1&0&-1&2&-1&0&-2&1 \\ 1&0&0&2&-1&-1&-2&1 \\
        2&-1&-1&4&-2&-1&-3&2 \\ 2&-1&-1&3&-1&-1&-2&1 \\
        1&-1&0&2&-1&0&-2&1 \\ 1&0&0&1&-1&0&-1&1 \\ 1&0&0&0&0&0&-1&1 \end{array}\right)$,
 $\left(\begin{array}{cccccccc} 1&1&-1&0&0&1&-2&1 \\ 0&1&0&0&0&1&-3&2 \\ 1&1&-1&0&0&2&-4&3 \\
        1&2&-1&0&0&2&-6&5 \\ 1&2&-1&0&0&1&-4&4 \\ 0&2&-1&0&0&1&-3&3 \\
        0&2&0&-1&0&1&-2&2 \\ 0&1&0&-1&1&0&-1&1 \end{array}\right)$,
  $\left(\begin{array}{cccccccc} 1&1&-1&0&1&-1&0&-1 \\ 1&1&-1&0&1&0&-1&-1 \\ 1&2&-2&0&2&-1&-1&-1 \\
        1&2&-2&0&3&-1&-2&-2 \\ 1&2&-1&-1&3&-1&-2&-1 \\
        0&2&0&-1&2&-1&-1&-1 \\ 0&1&0&-1&2&-1&0&-1 \\ 0&0&0&0&1&-1&0&0 \end{array}\right)$,
  $\left(\begin{array}{cccccccc} 1&1&0&0&-1&0&0&-1 \\ 2&1&-1&1&-2&1&-1&-1 \\
        2&2&-1&1&-2&0&-1&-1 \\ 3&2&-1&2&-4&1&-2&-1 \\
        3&1&-1&2&-3&0&-1&-1 \\ 2&1&0&1&-2&0&-1&-1 \\ 2&1&0&0&-1&0&-1&0 \\
        1&0&0&0&0&0&-1&0 \end{array}\right)$,
  $\left(\begin{array}{cccccccc} 2&0&-1&0&0&1&-1&-1 \\ 2&1&-1&0&0&1&-2&-1 \\ 3&0&-1&0&0&2&-3&-1 \\
        4&1&-2&0&0&3&-5&-1 \\ 4&1&-2&0&0&2&-4&0 \\ 3&1&-2&0&0&2&-3&0 \\
        2&1&-1&-1&1&1&-2&0 \\ 1&0&-1&0&1&0&-1&0 \end{array}\right)$.

The character table of $G^{s_{65}}$:\\
\noindent \begin{tabular}{c|ccccccccccccccccccccccc}
  & & & & & & & & & &10 & & & & & & & & & &20 & & &\\\hline
$\chi_{65}^{(1)}$&1&1&1&1&1&1&1&1&1&1&1&1&1&1&1 &1&1&1&1&1&1&1&1

\\$\chi_{65}^{(2)}$&1&-1&1&-1&1&-1&1&-1&-1&1&-1&1&-1&1&-1&-1
&1&-1&-1&1&1&-1&1

\\$\chi_{65}^{(3)}$&1&-1&1&-1&-1&1&1&-1&1&-1&1&-1&1&-1&1&1
&-1&1&1&1&1&1&1

\\$\chi_{65}^{(4)}$&1&1&1&1&-1&-1&1&1&-1&-1&-1&-1&-1&-1&-1&-1
&-1&-1&-1&1&1&-1&1

\\$\chi_{65}^{(5)}$&1&A&1&-1&-A&A&-/A&/A&A&1&-1&-/A&/A&1&-1&-1
&-/A&/A&/A&-A&1&-1&-A

\\$\chi_{65}^{(6)}$&1&/A&1&-1&-/A&/A&-A&A&/A&1&-1&-A&A&1&-1&-1
&-A&A&A&-/A&1&-1&-/A

\\$\chi_{65}^{(7)}$&1&A&1&-1&A&-A&-/A&/A&-A&-1&1&/A&-/A&-1&1&1
&/A&-/A&-/A&-A&1&1&-A

\\
\end {tabular}

\noindent \begin{tabular}{c|ccccccccccccccccccccccc}
  & & & & & & & & & &10 & & & & & & & & & &20 & & &\\\hline

\\$\chi_{65}^{({46})}$&2&.&A&.&.&/F&/B&.&-B&.&/E&.&H&.&-2&-/A
&.&G&-/C&-/G&-E&1&C

\\$\chi_{65}^{({47})}$&2&.&-1&.&.&-E&/A&.&A&.&2&.&-/E&.&2&-1
&.&-/E&/A&-E&2&-1&A

\\$\chi_{65}^{({48})}$&2&.&-1&.&.&-/E&A&.&/A&.&2&.&-E&.&2&-1
&.&-E&A&-/E&2&-1&/A

\\$\chi_{65}^{({49})}$&2&.&/A&.&.&-F&B&.&/B&.&-E&.&-/H&.&2&A
&.&-/G&C&-G&-/E&-1&/C

\\$\chi_{65}^{({50})}$&2&.&/A&.&.&-G&D&.&/D&.&-E&.&-/F&.&2&A
&.&-/H&B&-H&-/E&-1&/B

\\$\chi_{65}^{({51})}$&2&.&/A&.&.&-H&C&.&/C&.&-E&.&-/G&.&2&A
&.&-/F&D&-F&-/E&-1&/D

\\$\chi_{65}^{({52})}$&2&.&A&.&.&-/H&/C&.&C&.&-/E&.&-G&.&2&/A
&.&-F&/D&-/F&-E&-1&D

\\$\chi_{65}^{({53})}$&2&.&A&.&.&-/G&/D&.&D&.&-/E&.&-F&.&2&/A
&.&-H&/B&-/H&-E&-1&B

\\$\chi_{65}^{({54})}$&2&.&A&.&.&-/F&/B&.&B&.&-/E&.&-H&.&2&/A
&.&-G&/C&-/G&-E&-1&C\\
\end{tabular}
${}_{
 \!\!\!\!_{
\!\!\!\!_{ \!\!\!\!_{  \small
 \noindent
\begin{tabular}{c|ccccccccccccccccccccccc}
  & & & & & & &30 & & & & & & & & & & 40& & & & & &\\\hline
$\chi_{65}^{(1)}$&1&1&1&1&1&1&1&1&1&1&1&1&1&1&1 &1&1&1&1&1&1&1&1

\\$\chi_{65}^{(2)}$&-1&-1&1&1&-1&1&-1&-1&1&-1&1&-1&1&-1&1
&1&-1&1&-1&1&-1&-1&1

\\$\chi_{65}^{(3)}$&-1&1&1&-1&1&1&-1&1&1&-1&1&-1&1&-1&1
&-1&1&1&-1&-1&1&1&-1

\\$\chi_{65}^{(4)}$&1&-1&1&-1&-1&1&1&-1&1&1&1&1&1&1&1
&-1&-1&1&1&-1&-1&-1&-1

\\$\chi_{65}^{(5)}$&A&/A&-/A&-A&A&-/A&/A&A&1&-1&-/A&/A&1&-1&-A
&-/A&/A&-A&A&1&-1&/A&-A

\\$\chi_{65}^{(6)}$&/A&A&-A&-/A&/A&-A&A&/A&1&-1&-A&A&1&-1&-/A
&-A&A&-/A&/A&1&-1&A&-/A

\\$\chi_{65}^{(7)}$&A&-/A&-/A&A&-A&-/A&/A&-A&1&-1&-/A&/A&1&-1&-A
&/A&-/A&-A&A&-1&1&-/A&A

\\
\end{tabular}
 }}}}$

\noindent \begin{tabular}{c|ccccccccccccccccccccccc}
  & & & & & & &30 & & & & & & & & & & 40& & & & & &\\\hline
\\$\chi_{65}^{({36})}$&-C&-/D&-/B&-C&-C&-/D&-/D&-C&1&1&-/C&-/C&-/A&-/A&-C
&-/B&-/B&-B&-B&-A&-A&-/B&-D

\\$\chi_{65}^{({37})}$&.&1&2&.&-2&-1&.&1&-1&.&-1&.&-1&.&2
&.&-2&-1&.&.&-2&1&.

\\$\chi_{65}^{({38})}$&.&-1&2&.&2&-1&.&-1&-1&.&-1&.&-1&.&2
&.&2&-1&.&.&2&-1&.

\\$\chi_{65}^{({39})}$&.&-/A&-/E&.&E&/A&.&-A&-1&.&/A&.&-1&.&-E
&.&/E&A&.&.&-2&-/A&.

\\$\chi_{65}^{({40})}$&.&-A&-E&.&/E&A&.&-/A&-1&.&A&.&-1&.&-/E
&.&E&/A&.&.&-2&-A&.

\\$\chi_{65}^{({41})}$&.&-D&-/F&.&H&D&.&-/C&-1&.&C&.&A&.&-H
&.&/F&/B&.&.&/E&-B&.

\\$\chi_{65}^{({42})}$&.&-C&-/G&.&F&C&.&-/B&-1&.&B&.&A&.&-F
&.&/G&/D&.&.&/E&-D&.

\\$\chi_{65}^{({43})}$&.&-B&-/H&.&G&B&.&-/D&-1&.&D&.&A&.&-G
&.&/H&/C&.&.&/E&-C&.

\\$\chi_{65}^{({44})}$&.&-/B&-H&.&/G&/B&.&-D&-1&.&/D&.&/A&.&-/G
&.&H&C&.&.&E&-/C&.

\\$\chi_{65}^{({45})}$&.&-/C&-G&.&/F&/C&.&-B&-1&.&/B&.&/A&.&-/F
&.&G&D&.&.&E&-/D&.

\\$\chi_{65}^{({46})}$&.&-/D&-F&.&/H&/D&.&-C&-1&.&/C&.&/A&.&-/H
&.&F&B&.&.&E&-/B&.

\\$\chi_{65}^{({47})}$&.&/A&-/E&.&-E&/A&.&A&-1&.&/A&.&-1&.&-E
&.&-/E&A&.&.&2&/A&.

\\$\chi_{65}^{({48})}$&.&A&-E&.&-/E&A&.&/A&-1&.&A&.&-1&.&-/E
&.&-E&/A&.&.&2&A&.

\\$\chi_{65}^{({49})}$&.&D&-/F&.&-H&D&.&/C&-1&.&C&.&A&.&-H
&.&-/F&/B&.&.&-/E&B&.

\\$\chi_{65}^{({50})}$&.&C&-/G&.&-F&C&.&/B&-1&.&B&.&A&.&-F
&.&-/G&/D&.&.&-/E&D&.

\\$\chi_{65}^{({51})}$&.&B&-/H&.&-G&B&.&/D&-1&.&D&.&A&.&-G
&.&-/H&/C&.&.&-/E&C&.

\\$\chi_{65}^{({52})}$&.&/B&-H&.&-/G&/B&.&D&-1&.&/D&.&/A&.&-/G
&.&-H&C&.&.&-E&/C&.

\\$\chi_{65}^{({53})}$&.&/C&-G&.&-/F&/C&.&B&-1&.&/B&.&/A&.&-/F
&.&-G&D&.&.&-E&/D&.

\\$\chi_{65}^{({54})}$&.&/D&-F&.&-/H&/D&.&C&-1&.&/C&.&/A&.&-/H
&.&-F&B&.&.&-E&/B&.\\
\end{tabular}

\noindent \begin{tabular}{c|cccccccc}
  & & & & 50& & & &\\\hline
$\chi_{65}^{(1)}$&1&1&1&1&1&1&1&1

\\$\chi_{65}^{(2)}$&-1&1&1&-1&-1&1&1&1

\\$\chi_{65}^{(3)}$&1&1&1&1&1&1&1&1

\\$\chi_{65}^{(4)}$&-1&1&1&-1&-1&1&1&1

\\$\chi_{65}^{(5)}$&A&-/A&-A&-1&A&-/A&1&-A

\\$\chi_{65}^{(6)}$&/A&-A&-/A&-1&/A&-A&1&-/A

\\$\chi_{65}^{(7)}$&-A&-/A&-A&1&-A&-/A&1&-A

\\$\chi_{65}^{(8)}$&-/A&-A&-/A&1&-/A&-A&1&-/A

\\$\chi_{65}^{(9)}$&B&-/B&-B&A&B&-/D&-/A&-C

\\$\chi_{65}^{({10})}$&C&-/C&-C&A&C&-/B&-/A&-D

\\$\chi_{65}^{({11})}$&D&-/D&-D&A&D&-/C&-/A&-B

\\$\chi_{65}^{({12})}$&/B&-B&-/B&/A&/B&-D&-A&-/C

\\$\chi_{65}^{({13})}$&/D&-D&-/D&/A&/D&-C&-A&-/B

\\$\chi_{65}^{({14})}$&/C&-C&-/C&/A&/C&-B&-A&-/D

\\$\chi_{65}^{({15})}$&-/B&-B&-/B&-/A&-/B&-D&-A&-/C

\\$\chi_{65}^{({16})}$&-/D&-D&-/D&-/A&-/D&-C&-A&-/B

\\$\chi_{65}^{({17})}$&-/C&-C&-/C&-/A&-/C&-B&-A&-/D

\\$\chi_{65}^{({18})}$&-B&-/B&-B&-A&-B&-/D&-/A&-C

\\$\chi_{65}^{({19})}$&-C&-/C&-C&-A&-C&-/B&-/A&-D

\\$\chi_{65}^{({20})}$&-D&-/D&-D&-A&-D&-/C&-/A&-B

\\$\chi_{65}^{({21})}$&A&-/A&-A&-1&A&-/A&1&-A

\\$\chi_{65}^{({22})}$&/A&-A&-/A&-1&/A&-A&1&-/A

\\$\chi_{65}^{({23})}$&-A&-/A&-A&1&-A&-/A&1&-A

\\$\chi_{65}^{({24})}$&-/A&-A&-/A&1&-/A&-A&1&-/A

\\$\chi_{65}^{({25})}$&B&-/B&-B&A&B&-/D&-/A&-C
\end{tabular}

\noindent \begin{tabular}{c|cccccccc}
  & & & & 50& & & &\\\hline

\\$\chi_{65}^{({26})}$&C&-/C&-C&A&C&-/B&-/A&-D

\\$\chi_{65}^{({27})}$&D&-/D&-D&A&D&-/C&-/A&-B

\\$\chi_{65}^{({28})}$&/B&-B&-/B&/A&/B&-D&-A&-/C

\\$\chi_{65}^{({29})}$&/D&-D&-/D&/A&/D&-C&-A&-/B

\\$\chi_{65}^{({30})}$&/C&-C&-/C&/A&/C&-B&-A&-/D

\\$\chi_{65}^{({31})}$&-/B&-B&-/B&-/A&-/B&-D&-A&-/C

\\$\chi_{65}^{({32})}$&-/D&-D&-/D&-/A&-/D&-C&-A&-/B

\\$\chi_{65}^{({33})}$&-/C&-C&-/C&-/A&-/C&-B&-A&-/D

\\$\chi_{65}^{({34})}$&-B&-/B&-B&-A&-B&-/D&-/A&-C

\\$\chi_{65}^{({35})}$&-C&-/C&-C&-A&-C&-/B&-/A&-D

\\$\chi_{65}^{({36})}$&-D&-/D&-D&-A&-D&-/C&-/A&-B

\\$\chi_{65}^{({37})}$&-2&2&-1&1&1&2&2&2

\\$\chi_{65}^{({38})}$&2&2&-1&-1&-1&2&2&2

\\$\chi_{65}^{({39})}$&E&-/E&A&1&-A&-/E&2&-E

\\$\chi_{65}^{({40})}$&/E&-E&/A&1&-/A&-E&2&-/E

\\$\chi_{65}^{({41})}$&G&-/G&/D&-/A&-/D&-/H&-E&-F

\\$\chi_{65}^{({42})}$&H&-/H&/C&-/A&-/C&-/F&-E&-G

\\$\chi_{65}^{({43})}$&F&-/F&/B&-/A&-/B&-/G&-E&-H

\\$\chi_{65}^{({44})}$&/F&-F&B&-A&-B&-G&-/E&-/H

\\$\chi_{65}^{({45})}$&/H&-H&C&-A&-C&-F&-/E&-/G

\\$\chi_{65}^{({46})}$&/G&-G&D&-A&-D&-H&-/E&-/F

\\$\chi_{65}^{({47})}$&-E&-/E&A&-1&A&-/E&2&-E

\\$\chi_{65}^{({48})}$&-/E&-E&/A&-1&/A&-E&2&-/E

\\$\chi_{65}^{({49})}$&-G&-/G&/D&/A&/D&-/H&-E&-F

\\$\chi_{65}^{({50})}$&-H&-/H&/C&/A&/C&-/F&-E&-G

\\$\chi_{65}^{({51})}$&-F&-/F&/B&/A&/B&-/G&-E&-H

\\$\chi_{65}^{({52})}$&-/F&-F&B&A&B&-G&-/E&-/H

\\$\chi_{65}^{({53})}$&-/H&-H&C&A&C&-F&-/E&-/G

\\$\chi_{65}^{({54})}$&-/G&-G&D&A&D&-H&-/E&-/F\\
\end{tabular}

\noindent \noindent where A = -E(3)
  = (1-ER(-3))/2 = -b3,
B = -E(9)$^4$,C = -E(9)$^7$,D = E(9)$^4$+E(9)$^7$,E = -2*E(3)
  = 1-ER(-3) = 1-i3,
F = -2*E(9)$^5$,G = 2*E(9)$^2$+2*E(9)$^5$,H = -2*E(9)$^2$.

The generators of $G^{s_{66}}$ are:\\
$\left(\begin{array}{cccccccc} -1&1&2&-2& 0&1&0&0 \\
-2&2&2&-2&-1&2&0&0 \\ -3&2&4&-3&-1&2&
0&0 \\
        -4&3&5&-4&-2&4&-1&1 \\ -3&2&4&-3&-2&3&0&1 \\
        -2&1&3&-2&-2&3&0&0 \\ -1&1&2&-2&-1&2&0&0 \\ -1&0&1&-1&0&1&0&0 \end{array}\right)$,
 $\left(\begin{array}{cccccccc} 0&-1&0&1&0&0&-1&1 \\ -1&0&0&1&0&0&-1&1 \\ -1&-1&1&1&0&0&-1&1 \\
        -1&-1&0&2&0&0&-1&1 \\ -1&-1&0&1&1&0&-1&1 \\ -1&-1&0&1&0&1&-1&1 \\
        -1&-1&0&1&0&0&0&1 \\ 0&0&0&0&0&0&0&1 \end{array}\right)$,
  $\left(\begin{array}{cccccccc} 1&-1&0&1&0&-1&-1&1 \\ 2&-1&-1&1&1&-2&-1&2 \\ 2&-1&0&1&1&-3&-1&2 \\
        4&-2&-1&2&1&-4&-1&3 \\ 3&-2&-1&2&1&-3&-1&2 \\
        2&-2&0&1&1&-2&-1&2 \\ 1&-1&0&1&0&-1&-1&2 \\ 1&0&0&0&0&0&-1&1 \end{array}\right)$,
  $\left(\begin{array}{cccccccc} 1&-1&0&0&1&-1&0&1 \\ 0&-1&0&0&2&-2&0&2 \\ 0&-2&1&0&2&-2&0&2 \\
        0&-3&0&1&3&-3&0&3 \\ 0&-2&0&0&3&-2&0&2 \\ 0&-2&0&0&2&-1&0&2 \\
        0&-1&0&0&1&-1&1&1 \\ 0&0&0&0&0&0&0&1 \end{array}\right)$.

The representatives of conjugacy classes of   $G^{s_{66}}$ are:\\
$\left(\begin{array}{cccccccc} 1&0&0&0&0&0&0&0 \\ 0&1&0&0&
0&0&0&0 \\ 0&0&1&0&0&0&0&0 \\
        0&0&0&1&0&0&0&0 \\ 0&0&0&0&1&0&0&0 \\ 0&0&0&0&0&1&0&0 \\
        0&0&0&0&0&0&1&0 \\ 0&0&0&0&0&0&0&1 \end{array}\right)$,
  $\left(\begin{array}{cccccccc} -2&-1&1&1&-1&0&1&-1 \\ -2&-1&0&2&-2&1&1&-2 \\
        -3&-1&1&2&-2&0&2&-3 \\ -4&-2&1&3&-3&1&3&-5 \\
        -3&-2&1&2&-2&1&2&-4 \\ -2&-2&1&2&-2&1&1&-3 \\
        -1&-2&0&2&-1&0&1&-2 \\ 0&-1&0&1&-1&0&1&-1 \end{array}\right)$,
  $\left(\begin{array}{cccccccc} -2&0&0&0&0&1&0&0 \\ -3&-1&1&0&0&1&0&0 \\ -3&-1&0&0&0&2&0&0 \\
        -5&-1&1&-1&1&2&0&0 \\ -4&0&1&-1&0&2&0&0 \\ -3&0&0&0&0&1&0&0 \\
        -2&0&0&0&0&1&0&-1 \\ -1&0&0&0&0&0&1&-1 \end{array}\right)$,
  $\left(\begin{array}{cccccccc} -2&0&0&0&0&1&0&0 \\ -3&0&1&0&-1&2&0&-1 \\ -3&0&0&0&-1&3&0&-1 \\
        -5&1&1&-1&-1&4&0&-2 \\ -4&1&1&-1&-1&3&0&-1 \\
        -3&1&0&0&-1&2&0&-1 \\ -2&0&0&0&0&1&0&-1 \\ -1&0&0&0&0&0&1&-1 \end{array}\right)$,
  $\left(\begin{array}{cccccccc} -2&0&1&0&0&-1&1&1 \\ -3&-1&1&1&-1&0&1&1 \\ -4&0&1&1&-1&-1&2&1 \\
        -6&-1&2&2&-2&-1&3&1 \\ -5&-1&2&1&-1&-1&3&0 \\
        -4&-1&2&1&-1&-1&2&0 \\ -2&-1&1&1&-1&-1&2&0 \\
        -1&0&1&0&-1&0&1&0 \end{array}\right)$,
  $\left(\begin{array}{cccccccc} -2&0&1&0&0&-1&1&1 \\ -2&-1&1&0&0&-1&2&1 \\ -3&0&1&0&0&-2&3&1 \\
        -4&-1&2&0&0&-3&5&1 \\ -4&-1&2&0&0&-2&4&0 \\ -3&-1&2&0&0&-2&3&0 \\
        -2&-1&1&1&-1&-1&2&0 \\ -1&0&1&0&-1&0&1&0 \end{array}\right)$,
  $\left(\begin{array}{cccccccc} -2&0&2&-1&0&1&-1&1 \\ -3&0&2&-1&0&1&-1&2 \\ -4&0&4&-2&0&1&-1&2 \\
        -5&0&5&-3&0&2&-2&4 \\ -4&0&4&-2&-1&2&-1&3 \\
        -3&-1&3&-1&-1&2&-1&2 \\ -2&0&2&-1&-1&2&-1&1 \\
        -1&0&1&-1&0&1&0&0 \end{array}\right)$,
  $\left(\begin{array}{cccccccc} -2&0&2&-1&0&1&-1&1 \\ -3&1&2&-1&-1&2&-1&1 \\ -4&1&4&-2&-1&2&-1&1 \\
        -5&2&5&-3&-2&4&-2&2 \\ -4&1&4&-2&-2&3&-1&2 \\
        -3&0&3&-1&-2&3&-1&1 \\ -2&0&2&-1&-1&2&-1&1 \\
        -1&0&1&-1&0&1&0&0 \end{array}\right)$,
  $\left(\begin{array}{cccccccc} -2&1&1&0&-1&0&1&0 \\ -3&0&1&1&-2&1&1&0 \\ -4&1&1&1&-2&0&2&0 \\
        -6&0&2&2&-3&0&3&0 \\ -5&0&2&1&-2&0&3&-1 \\ -4&0&2&1&-2&0&2&-1 \\
        -2&0&1&1&-2&0&2&-1 \\ -1&0&1&0&-1&0&1&0 \end{array}\right)$,
  $\left(\begin{array}{cccccccc} -1&-1&0&0&1&0&0&1 \\ -3&-1&1&0&1&0&0&1 \\ -3&-2&1&0&1&1&0&1 \\
        -5&-2&1&0&2&1&0&1 \\ -4&-1&1&-1&2&1&0&1 \\ -3&-1&0&0&1&1&0&1 \\
        -2&-1&0&0&1&0&1&0 \\ -1&0&0&0&0&0&1&0 \end{array}\right)$,
  $\left(\begin{array}{cccccccc} -1&-1&0&0&1&0&0&1 \\ -2&-1&1&-1&2&-1&1&1 \\ -2&-2&1&-1&2&0&1&1 \\
        -3&-2&1&-2&4&-1&2&1 \\ -3&-1&1&-2&3&0&1&1 \\
        -2&-1&0&-1&2&0&1&1 \\ -2&-1&0&0&1&0&1&0 \\ -1&0&0&0&0&0&1&0 \end{array}\right)$,
   $\left(\begin{array}{cccccccc} -1&-1&1&0&-1&1&0&1 \\ -2&-1&1&1&-2&1&0&1 \\
        -2&-2&2&1&-3&2&0&1 \\ -3&-2&2&2&-5&3&0&2 \\ -2&-2&1&2&-4&2&1&1 \\
        -1&-2&0&2&-3&2&0&1 \\ 0&-1&0&1&-2&1&0&1 \\ 0&0&0&0&-1&1&0&0 \end{array}\right)$,
  $\left(\begin{array}{cccccccc} -1&-1&1&0&-1&1&0&1 \\ -1&-1&1&0&-1&0&1&1 \\
        -1&-2&2&0&-2&1&1&1 \\ -1&-2&2&0&-3&1&2&2 \\ -1&-2&1&1&-3&1&2&1 \\
        0&-2&0&1&-2&1&1&1 \\ 0&-1&0&1&-2&1&0&1 \\ 0&0&0&0&-1&1&0&0 \end{array}\right)$,
  $\left(\begin{array}{cccccccc} -1&0&0&0&0&0&0&0 \\ -1&-1&0&1&-1&1&-1&0 \\ -1&0&-1&1&-1&1&-1&0 \\
        -2&0&0&1&-2&2&-2&0 \\ -1&0&0&1&-2&1&-1&0 \\ -1&0&0&1&-1&0&-1&0 \\
        0&0&0&0&0&0&-1&0 \\ 0&0&0&0&0&0&0&-1 \end{array}\right)$,
  $\left(\begin{array}{cccccccc} -1&0&0&0&0&0&0&0 \\ 0&-1&0&0&0&0&0&0 \\ 0&0&-1&0&0&0&0&0 \\
        0&0&0&-1&0&0&0&0 \\ 0&0&0&0&-1&0&0&0 \\ 0&0&0&0&0&-1&0&0 \\
        0&0&0&0&0&0&-1&0 \\ 0&0&0&0&0&0&0&-1 \end{array}\right)$,
  $\left(\begin{array}{cccccccc} -1&0&0&0&0&1&0&0 \\ -3&0&1&0&0&1&0&0 \\ -3&-1&1&0&0&2&0&0 \\
        -5&-1&1&0&1&2&0&0 \\ -4&0&1&-1&1&2&0&0 \\ -3&0&0&0&0&2&0&0 \\
        -2&0&0&0&0&1&1&-1 \\ -1&0&0&0&0&0&1&0 \end{array}\right)$,
  $\left(\begin{array}{cccccccc} -1&0&1&-1&2&-1&-1&1 \\ -2&0&1&-1&2&0&-2&2 \\ -3&1&2&-2&3&-1&-2&2 \\
        -4&1&3&-3&4&0&-4&3 \\ -3&1&3&-3&3&0&-3&3 \\ -3&0&3&-2&2&0&-2&2 \\
        -2&0&2&-2&2&0&-1&1 \\ -1&0&1&-1&1&0&0&0 \end{array}\right)$,
  $\left(\begin{array}{cccccccc} -1&0&1&-1&2&-1&-1&1 \\ -1&0&1&-2&3&-1&-1&2 \\
        -2&1&2&-3&4&-2&-1&2 \\ -2&1&3&-5&6&-2&-2&3 \\
        -2&1&3&-4&4&-1&-2&3 \\ -2&0&3&-3&3&-1&-1&2 \\
        -2&0&2&-2&2&0&-1&1 \\ -1&0&1&-1&1&0&0&0 \end{array}\right)$,
  $\left(\begin{array}{cccccccc} -1&0&1&0&-2&2&0&0 \\ -2&0&1&1&-3&2&0&0 \\ -2&-1&2&1&-4&3&0&0 \\
        -3&-1&2&2&-6&4&0&1 \\ -2&-1&1&2&-5&3&1&0 \\ -1&-1&0&2&-4&3&0&0 \\
        0&0&0&1&-3&2&0&0 \\ 0&0&0&0&-1&1&0&0 \end{array}\right)$,
  $\left(\begin{array}{cccccccc} -1&0&1&0&-1&0&2&-2 \\ -1&0&0&1&-2&1&2&-3 \\ -2&0&1&1&-2&0&3&-4 \\
        -3&-1&1&2&-3&1&4&-6 \\ -2&-1&1&1&-2&1&3&-5 \\
        -1&-1&1&1&-2&1&2&-4 \\ 0&-1&0&1&-1&0&2&-3 \\
        0&-1&0&1&-1&0&1&-1 \end{array}\right)$,
  $\left(\begin{array}{cccccccc} -1&1&0&-1&0&1&1&-1 \\ -2&1&1&-1&-1&2&1&-2 \\ -2&1&0&-1&-1&3&1&-2 \\
        -4&2&1&-2&-1&4&1&-3 \\ -3&2&1&-2&-1&3&1&-2 \\
        -2&2&0&-1&-1&2&1&-2 \\ -1&1&0&-1&0&1&1&-2 \\ -1&0&0&0&0&0&1&-1 \end{array}\right)$,
 $\left(\begin{array}{cccccccc} -1&1&0&0&-1&1&0&-1 \\ -1&0&0&1&-2&2&-1&-1 \\
        -1&1&-1&1&-2&2&-1&-1 \\ -2&1&0&1&-3&3&-2&-1 \\
        -1&1&0&1&-3&2&-1&-1 \\ -1&1&0&1&-2&1&-1&-1 \\
        0&1&0&0&-1&1&-1&-1 \\ 0&0&0&0&0&0&0&-1 \end{array}\right)$,
  $\left(\begin{array}{cccccccc} -1&1&0&0&0&0&-1&2 \\ -2&0&1&0&0&0&-1&3 \\ -2&1&0&0&0&0&-1&4 \\
        -3&1&1&0&0&-1&-1&6 \\ -3&1&1&0&0&-1&0&4 \\ -3&1&1&0&0&-1&0&3 \\
        -2&1&1&0&-1&0&0&2 \\ -1&1&1&-1&0&0&0&1 \end{array}\right)$,
  $\left(\begin{array}{cccccccc} -1&1&0&0&0&0&-1&2 \\ -2&1&1&0&-1&1&-1&2 \\ -2&2&0&0&-1&1&-1&3 \\
        -3&3&1&0&-2&1&-1&4 \\ -3&2&1&0&-1&0&0&3 \\ -3&2&1&0&-1&0&0&2 \\
        -2&1&1&0&-1&0&0&2 \\ -1&1&1&-1&0&0&0&1 \end{array}\right)$,
  $\left(\begin{array}{cccccccc} -1&1&1&-1&1&0&-1&0 \\ -2&1&1&-1&1&1&-2&1 \\ -3&2&2&-2&2&0&-2&1 \\
        -4&2&3&-3&3&1&-4&2 \\ -3&2&3&-3&2&1&-3&2 \\ -3&1&3&-2&1&1&-2&1 \\
        -2&1&2&-2&1&1&-1&0 \\ -1&0&1&-1&1&0&0&0 \end{array}\right)$,
  $\left(\begin{array}{cccccccc} -1&1&2&-2&0&1&0&0 \\ -2&2&2&-2&-1&2&0&0 \\ -3&2&4&-3&-1&2&0&0 \\
        -4&3&5&-4&-2&4&-1&1 \\ -3&2&4&-3&-2&3&0&1 \\
        -2&1&3&-2&-2&3&0&0 \\ -1&1&2&-2&-1&2&0&0 \\ -1&0&1&-1&0&1&0&0 \end{array}\right)$,
  $\left(\begin{array}{cccccccc} 0&-2&0&1&0&0&0&-1 \\ 0&-2&-1&2&0&0&-1&-1 \\
        0&-3&0&2&0&0&-1&-2 \\ 0&-4&-1&3&0&1&-2&-3 \\
        1&-3&-1&2&0&1&-2&-2 \\ 1&-3&-1&2&0&1&-2&-1 \\
        1&-2&-1&1&1&0&-1&-1 \\ 1&-1&-1&1&0&0&0&-1 \end{array}\right)$,
  $\left(\begin{array}{cccccccc} 0&-2&0&1&0&0&0&-1 \\ 1&-2&-1&1&1&-1&0&-1 \\ 1&-3&0&1&1&-1&0&-2 \\
        2&-4&-1&1&2&-1&0&-3 \\ 2&-3&-1&1&1&0&-1&-2 \\
        2&-3&-1&1&1&0&-1&-1 \\ 1&-2&-1&1&1&0&-1&-1 \\
        1&-1&-1&1&0&0&0&-1 \end{array}\right)$,
  $\left(\begin{array}{cccccccc} 0&-1&-1&2&-2&1&0&0 \\ 0&-2&-1&3&-2&0&0&0 \\ 1&-3&-2&4&-3&1&0&0 \\
        1&-4&-3&6&-4&0&1&0 \\ 1&-3&-3&5&-3&0&1&-1 \\
        1&-2&-3&4&-2&0&0&0 \\ 1&-1&-2&3&-2&0&0&0 \\ 1&0&-1&1&-1&0&0&0 \end{array}\right)$,
  $\left(\begin{array}{cccccccc} 0&-1&-1&2&-2&1&0&0 \\ 0&-1&-1&3&-3&1&0&-1 \\
        1&-2&-2&4&-4&2&0&-1 \\ 1&-2&-3&6&-6&2&1&-2 \\
        1&-2&-3&5&-4&1&1&-2 \\ 1&-1&-3&4&-3&1&0&-1 \\
        1&-1&-2&3&-2&0&0&0 \\ 1&0&-1&1&-1&0&0&0 \end{array}\right)$,
  $\left(\begin{array}{cccccccc} 0&-1&0&1&-1&1&0&-2 \\ 0&-1&-1&2&-1&1&-1&-2 \\
        0&-2&0&2&-1&1&-1&-3 \\ 0&-3&-1&3&-1&2&-2&-4 \\
        1&-2&-1&2&-1&2&-2&-3 \\ 1&-2&-1&2&-1&2&-2&-2 \\
        1&-1&-1&1&0&1&-1&-2 \\ 1&-1&-1&1&0&0&0&-1 \end{array}\right)$,
  $\left(\begin{array}{cccccccc} 0&-1&0&1&0&0&-1&1 \\ -1&-1&0&1&1&-1&-1&2 \\ -1&-2&1&1&1&-1&-1&2 \\
        -1&-3&0&2&2&-2&-1&3 \\ -1&-2&0&1&2&-1&-1&2 \\
        -1&-2&0&1&1&0&-1&2 \\ -1&-1&0&1&0&0&0&1 \\ 0&0&0&0&0&0&0&1 \end{array}\right)$,
  $\left(\begin{array}{cccccccc} 0&-1&0&1&0&0&-1&1 \\ -1&0&0&1&0&0&-1&1 \\ -1&-1&1&1&0&0&-1&1 \\
        -1&-1&0&2&0&0&-1&1 \\ -1&-1&0&1&1&0&-1&1 \\ -1&-1&0&1&0&1&-1&1 \\
        -1&-1&0&1&0&0&0&1 \\ 0&0&0&0&0&0&0&1 \end{array}\right)$,
  $\left(\begin{array}{cccccccc} 0&0&-1&1&1&-1&-1&0 \\ 0&-1&-1&1&2&-1&-2&1 \\ 0&0&-2&1&3&-2&-2&1 \\
        0&-1&-2&1&5&-3&-3&1 \\ 0&0&-1&0&4&-2&-3&1 \\
        -1&0&0&0&3&-2&-2&1 \\ -1&0&0&0&2&-1&-1&0 \\ 0&0&0&0&1&-1&0&0 \end{array}\right)$,
  $\left(\begin{array}{cccccccc} 0&0&-1&1&1&-1&-1&0 \\ 0&0&-1&1&1&0&-2&0 \\ 0&1&-2&1&2&-1&-2&0 \\
        0&1&-2&1&3&-1&-3&-1 \\ 0&1&-1&0&3&-1&-3&0 \\
        -1&1&0&0&2&-1&-2&0 \\ -1&0&0&0&2&-1&-1&0 \\ 0&0&0&0&1&-1&0&0 \end{array}\right)$,
  $\left(\begin{array}{cccccccc} 0&0&0&1&-1&0&-1&0 \\ 1&-1&-1&2&-1&0&-2&1 \\
        1&0&-1&2&-1&-1&-2&1 \\ 2&-1&-1&3&-2&-1&-3&2 \\
        2&-1&-1&3&-2&-1&-2&1 \\ 1&-1&0&2&-1&-1&-2&1 \\
        1&0&0&1&-1&0&-2&1 \\ 1&0&0&0&0&0&-1&0 \end{array}\right)$,
  $\left(\begin{array}{cccccccc} 0&0&0&1&-1&0&-1&0 \\ 1&0&-1&2&-2&1&-2&0 \\ 1&1&-1&2&-2&0&-2&0 \\
        2&1&-1&3&-4&1&-3&0 \\ 2&0&-1&3&-3&0&-2&0 \\ 1&0&0&2&-2&0&-2&0 \\
        1&0&0&1&-1&0&-2&1 \\ 1&0&0&0&0&0&-1&0 \end{array}\right)$,
  $\left(\begin{array}{cccccccc} 0&2&0&-1&0&0&0&1 \\ -1&2&1&-1&-1&1&0&1 \\ -1&3&0&-1&-1&1&0&2 \\
        -2&4&1&-1&-2&1&0&3 \\ -2&3&1&-1&-1&0&1&2 \\ -2&3&1&-1&-1&0&1&1 \\
        -1&2&1&-1&-1&0&1&1 \\ -1&1&1&-1&0&0&0&1 \end{array}\right)$,
  $\left(\begin{array}{cccccccc} 1&-1&-2&2&0&-1&0&0 \\ 1&-2&-2&3&0&-1&-1&0 \\ 2&-2&-4&4&0&-1&-1&0 \\
        2&-3&-5&6&0&-2&-1&-1 \\ 2&-2&-4&4&1&-2&-1&-1 \\
        1&-1&-3&3&1&-2&-1&0 \\ 1&-1&-2&2&1&-2&0&0 \\ 1&0&-1&1&0&-1&0&0 \end{array}\right)$,
  $\left(\begin{array}{cccccccc} 1&-1&-2&2&0&-1&0&0 \\ 2&-2&-2&2&1&-2&0&0 \\
        3&-2&-4&3&1&-2&0&0 \\ 4&-3&-5&4&2&-4&1&-1 \\
        3&-2&-4&3&2&-3&0&-1 \\ 2&-1&-3&2&2&-3&0&0 \\
        1&-1&-2&2&1&-2&0&0 \\ 1&0&-1&1&0&-1&0&0 \end{array}\right)$,
  $\left(\begin{array}{cccccccc} 1&-1&-1&1&0&1&-2&0 \\ 1&-1&-1&1&1&0&-3&1 \\ 2&-2&-1&1&1&1&-4&1 \\
        3&-2&-2&1&2&1&-6&2 \\ 3&-1&-2&1&1&1&-5&2 \\ 2&-1&-2&1&1&1&-4&2 \\
        1&0&-1&0&1&1&-3&1 \\ 1&0&-1&0&1&0&-1&0 \end{array}\right)$,
  $\left(\begin{array}{cccccccc} 1&-1&-1&1&0&1&-2&0 \\ 1&0&-1&1&0&1&-3&0 \\ 2&-1&-1&1&0&2&-4&0 \\
        3&0&-2&1&0&3&-6&0 \\ 3&0&-2&1&0&2&-5&1 \\ 2&0&-2&1&0&2&-4&1 \\
        1&0&-1&0&1&1&-3&1 \\ 1&0&-1&0&1&0&-1&0 \end{array}\right)$,
  $\left(\begin{array}{cccccccc} 1&-1&0&1&0&-1&-1&1 \\ 1&-1&-1&2&0&-1&-2&2 \\ 1&-1&0&2&0&-2&-2&2 \\
        2&-2&-1&4&-1&-2&-3&3 \\ 2&-2&-1&3&0&-2&-2&2 \\
        1&-2&0&2&0&-1&-2&2 \\ 1&-1&0&1&0&-1&-1&2 \\ 1&0&0&0&0&0&-1&1 \end{array}\right)$,
 $\left(\begin{array}{cccccccc} 1&-1&0&1&0&-1&-1&1 \\ 2&-1&-1&1&1&-2&-1&2 \\
        2&-1&0&1&1&-3&-1&2 \\ 4&-2&-1&2&1&-4&-1&3 \\
        3&-2&-1&2&1&-3&-1&2 \\ 2&-2&0&1&1&-2&-1&2 \\
        1&-1&0&1&0&-1&-1&2 \\ 1&0&0&0&0&0&-1&1 \end{array}\right)$,
  $\left(\begin{array}{cccccccc} 1&0&-2&2&-1&0&0&-1 \\ 1&-1&-2&3&-1&0&-1&-1 \\
        2&-1&-4&4&-1&0&-1&-1 \\ 2&-2&-5&6&-1&-1&-1&-2 \\
        2&-1&-4&4&0&-1&-1&-2 \\ 1&0&-3&3&0&-1&-1&-1 \\
        1&0&-2&2&0&-1&0&-1 \\ 1&0&-1&1&0&-1&0&0 \end{array}\right)$,
  $\left(\begin{array}{cccccccc} 1&0&-1&0&1&0&-2&2 \\ 0&0&0&0&1&0&-3&3 \\ 1&0&-1&0&1&1&-4&4 \\
        1&1&-1&0&1&1&-6&6 \\ 1&1&-1&0&1&0&-4&5 \\ 0&1&-1&0&1&0&-3&4 \\
        0&1&0&-1&1&0&-2&3 \\ 0&1&0&-1&1&0&-1&1 \end{array}\right)$,
  $\left(\begin{array}{cccccccc} 1&0&-1&0&1&0&-2&2 \\ 1&0&0&-1&2&-1&-2&3 \\ 2&0&-1&-1&2&0&-3&4 \\
        3&1&-1&-2&3&-1&-4&6 \\ 2&1&-1&-1&2&-1&-3&5 \\
        1&1&-1&-1&2&-1&-2&4 \\ 0&1&0&-1&1&0&-2&3 \\ 0&1&0&-1&1&0&-1&1 \end{array}\right)$,
  $\left(\begin{array}{cccccccc} 1&0&-1&1&-2&1&1&-1 \\ 1&0&-1&2&-3&1&1&-2 \\
        2&-1&-2&3&-4&2&1&-2 \\ 2&-1&-3&5&-6&2&2&-3 \\
        2&-1&-3&4&-4&1&2&-3 \\ 2&0&-3&3&-3&1&1&-2 \\
        2&0&-2&2&-2&0&1&-1 \\ 1&0&-1&1&-1&0&0&0 \end{array}\right)$,
  $\left(\begin{array}{cccccccc} -2&-1&1&1&-1&0&1&-1 \\ -2&-2&0&2&-1&0&1&-1 \\
        -3&-2&1&2&-1&-1&2&-2 \\ -4&-4&1&3&-1&-1&3&-3 \\
        -3&-3&1&2&-1&0&2&-3 \\ -2&-3&1&2&-1&0&1&-2 \\
        -1&-2&0&2&-1&0&1&-2 \\ 0&-1&0&1&-1&0&1&-1 \end{array}\right)$,
  $\left(\begin{array}{cccccccc} 1&0&0&1&-1&0&-1&0 \\ 1&0&-1&2&-1&0&-2&1 \\ 1&0&0&2&-1&-1&-2&1 \\
        2&-1&-1&4&-2&-1&-3&2 \\ 2&-1&-1&3&-1&-1&-2&1 \\
        1&-1&0&2&-1&0&-2&1 \\ 1&0&0&1&-1&0&-1&1 \\ 1&0&0&0&0&0&-1&1 \end{array}\right)$,
 $\left(\begin{array}{cccccccc} 1&1&-1&0&0&1&-2&1 \\ 0&1&0&0&0&1&-3&2 \\ 1&1&-1&0&0&2&-4&3 \\
        1&2&-1&0&0&2&-6&5 \\ 1&2&-1&0&0&1&-4&4 \\ 0&2&-1&0&0&1&-3&3 \\
        0&2&0&-1&0&1&-2&2 \\ 0&1&0&-1&1&0&-1&1 \end{array}\right)$,
  $\left(\begin{array}{cccccccc} 1&1&-1&0&1&-1&0&-1 \\ 1&1&-1&0&1&0&-1&-1 \\ 1&2&-2&0&2&-1&-1&-1 \\
        1&2&-2&0&3&-1&-2&-2 \\ 1&2&-1&-1&3&-1&-2&-1 \\
        0&2&0&-1&2&-1&-1&-1 \\ 0&1&0&-1&2&-1&0&-1 \\ 0&0&0&0&1&-1&0&0 \end{array}\right)$,
  $\left(\begin{array}{cccccccc} 1&1&0&0&-1&0&0&-1 \\ 2&1&-1&1&-2&1&-1&-1 \\
        2&2&-1&1&-2&0&-1&-1 \\ 3&2&-1&2&-4&1&-2&-1 \\
        3&1&-1&2&-3&0&-1&-1 \\ 2&1&0&1&-2&0&-1&-1 \\ 2&1&0&0&-1&0&-1&0 \\
        1&0&0&0&0&0&-1&0 \end{array}\right)$,
  $\left(\begin{array}{cccccccc} 2&0&-1&0&0&1&-1&-1 \\ 2&1&-1&0&0&1&-2&-1 \\ 3&0&-1&0&0&2&-3&-1 \\
        4&1&-2&0&0&3&-5&-1 \\ 4&1&-2&0&0&2&-4&0 \\ 3&1&-2&0&0&2&-3&0 \\
        2&1&-1&-1&1&1&-2&0 \\ 1&0&-1&0&1&0&-1&0 \end{array}\right)$.

The character table of $G^{s_{66}}$:\\
\noindent \begin{tabular}{c|ccccccccccccccccccccccc}
  & & & & & & & & & &10 & & & & & & & & & &20 & & &\\\hline
$\chi_{66}^{(1)}$&1&1&1&1&1&1&1&1&1&1&1&1&1&1&1 &1&1&1&1&1&1&1&1

\\$\chi_{66}^{(2)}$&1&-1&1&-1&1&-1&1&-1&-1&1&-1&1&-1&1&-1&-1
&1&-1&-1&1&1&-1&1

\\$\chi_{66}^{(3)}$&1&-1&1&-1&-1&1&1&-1&1&-1&1&-1&1&-1&1&1
&-1&1&1&1&1&1&1

\\$\chi_{66}^{(4)}$&1&1&1&1&-1&-1&1&1&-1&-1&-1&-1&-1&-1&-1&-1
&-1&-1&-1&1&1&-1&1

\\$\chi_{66}^{(5)}$&1&A&1&-1&-A&A&-/A&/A&A&1&-1&-/A&/A&1&-1&-1
&-/A&/A&/A&-A&1&-1&-A

\\$\chi_{66}^{(6)}$&1&/A&1&-1&-/A&/A&-A&A&/A&1&-1&-A&A&1&-1&-1
&-A&A&A&-/A&1&-1&-/A

\\$\chi_{66}^{(7)}$&1&A&1&-1&A&-A&-/A&/A&-A&-1&1&/A&-/A&-1&1&1
&/A&-/A&-/A&-A&1&1&-A

\\
\end {tabular}

\noindent \begin{tabular}{c|ccccccccccccccccccccccc}
  & & & & & & & & & &10 & & & & & & & & & &20 & & &\\\hline

\\$\chi_{66}^{({46})}$&2&.&A&.&.&/F&/B&.&-B&.&/E&.&H&.&-2&-/A
&.&G&-/C&-/G&-E&1&C

\\$\chi_{66}^{({47})}$&2&.&-1&.&.&-E&/A&.&A&.&2&.&-/E&.&2&-1
&.&-/E&/A&-E&2&-1&A

\\$\chi_{66}^{({48})}$&2&.&-1&.&.&-/E&A&.&/A&.&2&.&-E&.&2&-1
&.&-E&A&-/E&2&-1&/A

\\$\chi_{66}^{({49})}$&2&.&/A&.&.&-F&B&.&/B&.&-E&.&-/H&.&2&A
&.&-/G&C&-G&-/E&-1&/C

\\$\chi_{66}^{({50})}$&2&.&/A&.&.&-G&D&.&/D&.&-E&.&-/F&.&2&A
&.&-/H&B&-H&-/E&-1&/B

\\$\chi_{66}^{({51})}$&2&.&/A&.&.&-H&C&.&/C&.&-E&.&-/G&.&2&A
&.&-/F&D&-F&-/E&-1&/D

\\$\chi_{66}^{({52})}$&2&.&A&.&.&-/H&/C&.&C&.&-/E&.&-G&.&2&/A
&.&-F&/D&-/F&-E&-1&D

\\$\chi_{66}^{({53})}$&2&.&A&.&.&-/G&/D&.&D&.&-/E&.&-F&.&2&/A
&.&-H&/B&-/H&-E&-1&B

\\$\chi_{66}^{({54})}$&2&.&A&.&.&-/F&/B&.&B&.&-/E&.&-H&.&2&/A
&.&-G&/C&-/G&-E&-1&C\\
\end{tabular}
${}_{
 \!\!\!\!_{
\!\!\!\!\!_{ \!\!\!\!_{  \small
 \noindent
\begin{tabular}{c|ccccccccccccccccccccccc}
  & & & & & & & 30& & & & & & & & & &40 & & & & & &\\\hline
$\chi_{66}^{(1)}$&1&1&1&1&1&1&1&1&1&1&1&1&1&1&1 &1&1&1&1&1&1&1&1

\\$\chi_{66}^{(2)}$&-1&-1&1&1&-1&1&-1&-1&1&-1&1&-1&1&-1&1
&1&-1&1&-1&1&-1&-1&1

\\$\chi_{66}^{(3)}$&-1&1&1&-1&1&1&-1&1&1&-1&1&-1&1&-1&1
&-1&1&1&-1&-1&1&1&-1

\\$\chi_{66}^{(4)}$&1&-1&1&-1&-1&1&1&-1&1&1&1&1&1&1&1
&-1&-1&1&1&-1&-1&-1&-1

\\$\chi_{66}^{(5)}$&A&/A&-/A&-A&A&-/A&/A&A&1&-1&-/A&/A&1&-1&-A
&-/A&/A&-A&A&1&-1&/A&-A

\\$\chi_{66}^{(6)}$&/A&A&-A&-/A&/A&-A&A&/A&1&-1&-A&A&1&-1&-/A
&-A&A&-/A&/A&1&-1&A&-/A

\\$\chi_{66}^{(7)}$&A&-/A&-/A&A&-A&-/A&/A&-A&1&-1&-/A&/A&1&-1&-A
&/A&-/A&-A&A&-1&1&-/A&A

\\
\end{tabular}
 }}}}$

\noindent \begin{tabular}{c|ccccccccccccccccccccccc}
  & & & & & & & 30& & & & & & & & & &40 & & & & & &\\\hline
\\$\chi_{66}^{({36})}$&-C&-/D&-/B&-C&-C&-/D&-/D&-C&1&1&-/C&-/C&-/A&-/A&-C
&-/B&-/B&-B&-B&-A&-A&-/B&-D

\\$\chi_{66}^{({37})}$&.&1&2&.&-2&-1&.&1&-1&.&-1&.&-1&.&2
&.&-2&-1&.&.&-2&1&.

\\$\chi_{66}^{({38})}$&.&-1&2&.&2&-1&.&-1&-1&.&-1&.&-1&.&2
&.&2&-1&.&.&2&-1&.

\\$\chi_{66}^{({39})}$&.&-/A&-/E&.&E&/A&.&-A&-1&.&/A&.&-1&.&-E
&.&/E&A&.&.&-2&-/A&.

\\$\chi_{66}^{({40})}$&.&-A&-E&.&/E&A&.&-/A&-1&.&A&.&-1&.&-/E
&.&E&/A&.&.&-2&-A&.

\\$\chi_{66}^{({41})}$&.&-D&-/F&.&H&D&.&-/C&-1&.&C&.&A&.&-H
&.&/F&/B&.&.&/E&-B&.

\\$\chi_{66}^{({42})}$&.&-C&-/G&.&F&C&.&-/B&-1&.&B&.&A&.&-F
&.&/G&/D&.&.&/E&-D&.

\\$\chi_{66}^{({43})}$&.&-B&-/H&.&G&B&.&-/D&-1&.&D&.&A&.&-G
&.&/H&/C&.&.&/E&-C&.

\\$\chi_{66}^{({44})}$&.&-/B&-H&.&/G&/B&.&-D&-1&.&/D&.&/A&.&-/G
&.&H&C&.&.&E&-/C&.

\\$\chi_{66}^{({45})}$&.&-/C&-G&.&/F&/C&.&-B&-1&.&/B&.&/A&.&-/F
&.&G&D&.&.&E&-/D&.

\\$\chi_{66}^{({46})}$&.&-/D&-F&.&/H&/D&.&-C&-1&.&/C&.&/A&.&-/H
&.&F&B&.&.&E&-/B&.

\\$\chi_{66}^{({47})}$&.&/A&-/E&.&-E&/A&.&A&-1&.&/A&.&-1&.&-E
&.&-/E&A&.&.&2&/A&.

\\$\chi_{66}^{({48})}$&.&A&-E&.&-/E&A&.&/A&-1&.&A&.&-1&.&-/E
&.&-E&/A&.&.&2&A&.

\\$\chi_{66}^{({49})}$&.&D&-/F&.&-H&D&.&/C&-1&.&C&.&A&.&-H
&.&-/F&/B&.&.&-/E&B&.

\\$\chi_{66}^{({50})}$&.&C&-/G&.&-F&C&.&/B&-1&.&B&.&A&.&-F
&.&-/G&/D&.&.&-/E&D&.

\\$\chi_{66}^{({51})}$&.&B&-/H&.&-G&B&.&/D&-1&.&D&.&A&.&-G
&.&-/H&/C&.&.&-/E&C&.

\\$\chi_{66}^{({52})}$&.&/B&-H&.&-/G&/B&.&D&-1&.&/D&.&/A&.&-/G
&.&-H&C&.&.&-E&/C&.

\\$\chi_{66}^{({53})}$&.&/C&-G&.&-/F&/C&.&B&-1&.&/B&.&/A&.&-/F
&.&-G&D&.&.&-E&/D&.

\\$\chi_{66}^{({54})}$&.&/D&-F&.&-/H&/D&.&C&-1&.&/C&.&/A&.&-/H
&.&-F&B&.&.&-E&/B&.\\
\end{tabular}

\noindent \begin{tabular}{c|cccccccc}
  & & & & 50& & & &\\\hline
$\chi_{66}^{(1)}$&1&1&1&1&1&1&1&1

\\$\chi_{66}^{(2)}$&-1&1&1&-1&-1&1&1&1

\\$\chi_{66}^{(3)}$&1&1&1&1&1&1&1&1

\\$\chi_{66}^{(4)}$&-1&1&1&-1&-1&1&1&1

\\$\chi_{66}^{(5)}$&A&-/A&-A&-1&A&-/A&1&-A

\\$\chi_{66}^{(6)}$&/A&-A&-/A&-1&/A&-A&1&-/A

\\$\chi_{66}^{(7)}$&-A&-/A&-A&1&-A&-/A&1&-A

\\$\chi_{66}^{(8)}$&-/A&-A&-/A&1&-/A&-A&1&-/A

\\$\chi_{66}^{(9)}$&B&-/B&-B&A&B&-/D&-/A&-C

\\$\chi_{66}^{({10})}$&C&-/C&-C&A&C&-/B&-/A&-D

\\$\chi_{66}^{({11})}$&D&-/D&-D&A&D&-/C&-/A&-B

\\$\chi_{66}^{({12})}$&/B&-B&-/B&/A&/B&-D&-A&-/C

\\$\chi_{66}^{({13})}$&/D&-D&-/D&/A&/D&-C&-A&-/B

\\$\chi_{66}^{({14})}$&/C&-C&-/C&/A&/C&-B&-A&-/D

\\$\chi_{66}^{({15})}$&-/B&-B&-/B&-/A&-/B&-D&-A&-/C

\\$\chi_{66}^{({16})}$&-/D&-D&-/D&-/A&-/D&-C&-A&-/B

\\$\chi_{66}^{({17})}$&-/C&-C&-/C&-/A&-/C&-B&-A&-/D

\\$\chi_{66}^{({18})}$&-B&-/B&-B&-A&-B&-/D&-/A&-C

\\$\chi_{66}^{({19})}$&-C&-/C&-C&-A&-C&-/B&-/A&-D

\\$\chi_{66}^{({20})}$&-D&-/D&-D&-A&-D&-/C&-/A&-B

\\$\chi_{66}^{({21})}$&A&-/A&-A&-1&A&-/A&1&-A

\\$\chi_{66}^{({22})}$&/A&-A&-/A&-1&/A&-A&1&-/A

\\$\chi_{66}^{({23})}$&-A&-/A&-A&1&-A&-/A&1&-A

\\$\chi_{66}^{({24})}$&-/A&-A&-/A&1&-/A&-A&1&-/A

\\$\chi_{66}^{({25})}$&B&-/B&-B&A&B&-/D&-/A&-C
\end{tabular}

\noindent \begin{tabular}{c|cccccccc}
  & & & & 50& & & &\\\hline

\\$\chi_{66}^{({26})}$&C&-/C&-C&A&C&-/B&-/A&-D

\\$\chi_{66}^{({27})}$&D&-/D&-D&A&D&-/C&-/A&-B

\\$\chi_{66}^{({28})}$&/B&-B&-/B&/A&/B&-D&-A&-/C

\\$\chi_{66}^{({29})}$&/D&-D&-/D&/A&/D&-C&-A&-/B

\\$\chi_{66}^{({30})}$&/C&-C&-/C&/A&/C&-B&-A&-/D

\\$\chi_{66}^{({31})}$&-/B&-B&-/B&-/A&-/B&-D&-A&-/C

\\$\chi_{66}^{({32})}$&-/D&-D&-/D&-/A&-/D&-C&-A&-/B

\\$\chi_{66}^{({33})}$&-/C&-C&-/C&-/A&-/C&-B&-A&-/D

\\$\chi_{66}^{({34})}$&-B&-/B&-B&-A&-B&-/D&-/A&-C

\\$\chi_{66}^{({35})}$&-C&-/C&-C&-A&-C&-/B&-/A&-D

\\$\chi_{66}^{({36})}$&-D&-/D&-D&-A&-D&-/C&-/A&-B

\\$\chi_{66}^{({37})}$&-2&2&-1&1&1&2&2&2

\\$\chi_{66}^{({38})}$&2&2&-1&-1&-1&2&2&2

\\$\chi_{66}^{({39})}$&E&-/E&A&1&-A&-/E&2&-E

\\$\chi_{66}^{({40})}$&/E&-E&/A&1&-/A&-E&2&-/E

\\$\chi_{66}^{({41})}$&G&-/G&/D&-/A&-/D&-/H&-E&-F

\\$\chi_{66}^{({42})}$&H&-/H&/C&-/A&-/C&-/F&-E&-G

\\$\chi_{66}^{({43})}$&F&-/F&/B&-/A&-/B&-/G&-E&-H

\\$\chi_{66}^{({44})}$&/F&-F&B&-A&-B&-G&-/E&-/H

\\$\chi_{66}^{({45})}$&/H&-H&C&-A&-C&-F&-/E&-/G

\\$\chi_{66}^{({46})}$&/G&-G&D&-A&-D&-H&-/E&-/F

\\$\chi_{66}^{({47})}$&-E&-/E&A&-1&A&-/E&2&-E

\\$\chi_{66}^{({48})}$&-/E&-E&/A&-1&/A&-E&2&-/E

\\$\chi_{66}^{({49})}$&-G&-/G&/D&/A&/D&-/H&-E&-F

\\$\chi_{66}^{({50})}$&-H&-/H&/C&/A&/C&-/F&-E&-G

\\$\chi_{66}^{({51})}$&-F&-/F&/B&/A&/B&-/G&-E&-H

\\$\chi_{66}^{({52})}$&-/F&-F&B&A&B&-G&-/E&-/H

\\$\chi_{66}^{({53})}$&-/H&-H&C&A&C&-F&-/E&-/G

\\$\chi_{66}^{({54})}$&-/G&-G&D&A&D&-H&-/E&-/F\\
\end{tabular}

\noindent \noindent where A = -E(3)
  = (1-ER(-3))/2 = -b3,
B = -E(9)$^4$ ,C = -E(9)$^7$,D = E(9)$^4$+E(9)$^7$,E = -2*E(3)
  = 1-ER(-3) = 1-i3,
F = -2*E(9)$^5$,G = 2*E(9)$^2$+2*E(9)$^5$,H = -2*E(9)$^2$.

The generators of $G^{s_{67}}$ are:\\
 $\left(\begin{array}{cccccccc}
0&-1&0&1& -1&1&0&-2 \\ 1&0&-1&1&-1&1&0&-3 \\ 1&-1&0&1&-1&1&
0&-4 \\
        2&-1&-1&1&-1&2&0&-6 \\ 2&-1&-1&1&-1&2&-1&-4 \\
        2&-1&-1&1&-1&2&-1&-3 \\ 1&-1&-1&1&0&1&-1&-2 \\
        1&-1&-1&1&0&0&0&-1 \end{array}\right)$,
  $\left(\begin{array}{cccccccc} 1&-1&0&1&0&-1&-1&1 \\ 2&-1&-1&1&1&-2&-1&2 \\ 2&-1&0&1&1&-3&-1&2 \\
        4&-2&-1&2&1&-4&-1&3 \\ 3&-2&-1&2&1&-3&-1&2 \\
        2&-2&0&1&1&-2&-1&2 \\ 1&-1&0&1&0&-1&-1&2 \\ 1&0&0&0&0&0&-1&1
        \end{array}\right)$.

The representatives of conjugacy classes of   $G^{s_{67}}$ are:\\
$\left(\begin{array}{cccccccc} 1&0&0&0&0&0&0&0 \\ 0&1&0&0&
0&0&0&0 \\ 0&0&1&0&0&0&0&0 \\
        0&0&0&1&0&0&0&0 \\ 0&0&0&0&1&0&0&0 \\ 0&0&0&0&0&1&0&0 \\
        0&0&0&0&0&0&1&0 \\ 0&0&0&0&0&0&0&1 \end{array}\right)$,
  $\left(\begin{array}{cccccccc} -2&1&1&0&-1&0&1&0 \\ -2&1&1&0&-2&1&2&-1 \\ -3&2&1&0&-2&0&3&-1 \\
        -4&2&2&0&-3&0&5&-2 \\ -4&1&2&0&-2&0&4&-2 \\ -3&1&2&0&-2&0&3&-2 \\
        -2&0&1&1&-2&0&2&-1 \\ -1&0&1&0&-1&0&1&0 \end{array}\right)$,
  $\left(\begin{array}{cccccccc} -1&-1&0&0&1&0&0&1 \\ -2&-1&1&-1&2&-1&1&1 \\ -2&-2&1&-1&2&0&1&1 \\
        -3&-2&1&-2&4&-1&2&1 \\ -3&-1&1&-2&3&0&1&1 \\
        -2&-1&0&-1&2&0&1&1 \\ -2&-1&0&0&1&0&1&0 \\ -1&0&0&0&0&0&1&0
        \end{array}\right)$,\\
   $\left(\begin{array}{cccccccc} -1&-1&1&0&-1&1&0&1 \\ -1&-1&1&0&-1&0&1&1 \\
        -1&-2&2&0&-2&1&1&1 \\ -1&-2&2&0&-3&1&2&2 \\ -1&-2&1&1&-3&1&2&1 \\
        0&-2&0&1&-2&1&1&1 \\ 0&-1&0&1&-2&1&0&1 \\ 0&0&0&0&-1&1&0&0 \end{array}\right)$,
  $\left(\begin{array}{cccccccc} -1&-1&1&0&0&-1&2&-1 \\ -1&-2&0&1&0&-1&2&-1 \\
        -2&-2&1&1&0&-2&3&-2 \\ -3&-4&1&2&0&-2&4&-3 \\
        -2&-3&1&1&0&-1&3&-3 \\ -1&-3&1&1&0&-1&2&-2 \\
        0&-2&0&1&0&-1&2&-2 \\ 0&-1&0&1&-1&0&1&-1 \end{array}\right)$,
  $\left(\begin{array}{cccccccc} -1&0&0&-1&1&0&1&0 \\ -2&-1&1&-1&1&0&1&0 \\ -2&-1&0&-1&1&1&1&0 \\
        -4&-1&1&-2&2&1&1&0 \\ -3&0&1&-2&1&1&1&0 \\ -2&0&0&-1&1&0&1&0 \\
        -1&0&0&-1&1&0&1&-1 \\ -1&0&0&0&0&0&1&-1
        \end{array}\right)$,\\
  $\left(\begin{array}{cccccccc} -1&0&0&0&0&0&0&0 \\ 0&-1&0&0&0&0&0&0 \\ 0&0&-1&0&0&0&0&0 \\
        0&0&0&-1&0&0&0&0 \\ 0&0&0&0&-1&0&0&0 \\ 0&0&0&0&0&-1&0&0 \\
        0&0&0&0&0&0&-1&0 \\ 0&0&0&0&0&0&0&-1 \end{array}\right)$,
  $\left(\begin{array}{cccccccc} -1&0&0&0&0&1&0&0 \\ -2&1&1&-1&0&1&1&-1 \\ -2&0&1&-1&0&2&1&-1 \\
        -3&1&1&-2&1&2&2&-2 \\ -3&1&1&-2&1&2&1&-1 \\ -2&1&0&-1&0&2&1&-1 \\
        -2&0&0&0&0&1&1&-1 \\ -1&0&0&0&0&0&1&0 \end{array}\right)$,
  $\left(\begin{array}{cccccccc} -1&0&1&-1&2&-1&-1&1 \\ -1&0&1&-2&3&-1&-1&2 \\
        -2&1&2&-3&4&-2&-1&2 \\ -2&1&3&-5&6&-2&-2&3 \\
        -2&1&3&-4&4&-1&-2&3 \\ -2&0&3&-3&3&-1&-1&2 \\
        -2&0&2&-2&2&0&-1&1 \\ -1&0&1&-1&1&0&0&0
        \end{array}\right)$,\\
  $\left(\begin{array}{cccccccc} -1&0&1&0&-2&2&0&0 \\ -1&1&1&0&-3&2&1&-1 \\ -1&0&2&0&-4&3&1&-1 \\
        -1&1&2&0&-6&4&2&-1 \\ -1&0&1&1&-5&3&2&-1 \\ 0&0&0&1&-4&3&1&-1 \\
        0&0&0&1&-3&2&0&0 \\ 0&0&0&0&-1&1&0&0 \end{array}\right)$,
  $\left(\begin{array}{cccccccc} -1&0&1&0&-1&0&2&-2 \\ -1&0&0&1&-2&1&2&-3 \\ -2&0&1&1&-2&0&3&-4 \\
        -3&-1&1&2&-3&1&4&-6 \\ -2&-1&1&1&-2&1&3&-5 \\
        -1&-1&1&1&-2&1&2&-4 \\ 0&-1&0&1&-1&0&2&-3 \\
        0&-1&0&1&-1&0&1&-1 \end{array}\right)$,
  $\left(\begin{array}{cccccccc} -1&0&2&-2&1&0&0&1 \\ -2&0&2&-2&1&0&0&2 \\ -3&0&4&-3&1&0&0&2 \\
        -4&0&5&-4&1&1&-1&4 \\ -3&0&4&-3&0&1&0&3 \\ -2&-1&3&-2&0&1&0&2 \\
        -1&0&2&-2&0&1&0&1 \\ -1&0&1&-1&0&1&0&0
        \end{array}\right)$,\\
  $\left(\begin{array}{cccccccc} -1&1&0&-1&0&1&1&-1 \\ -2&1&1&-1&-1&2&1&-2 \\ -2&1&0&-1&-1&3&1&-2 \\
        -4&2&1&-2&-1&4&1&-3 \\ -3&2&1&-2&-1&3&1&-2 \\
        -2&2&0&-1&-1&2&1&-2 \\ -1&1&0&-1&0&1&1&-2 \\ -1&0&0&0&0&0&1&-1 \end{array}\right)$,
  $\left(\begin{array}{cccccccc} -1&1&0&0&-1&1&0&-1 \\ 0&1&0&0&-2&2&0&-2 \\ 0&2&-1&0&-2&2&0&-2 \\
        0&3&0&-1&-3&3&0&-3 \\ 0&2&0&0&-3&2&0&-2 \\ 0&2&0&0&-2&1&0&-2 \\
        0&1&0&0&-1&1&-1&-1 \\ 0&0&0&0&0&0&0&-1 \end{array}\right)$,
  $\left(\begin{array}{cccccccc} -1&1&1&-1&1&0&-1&0 \\ -1&2&1&-2&1&1&-1&0 \\ -2&3&2&-3&2&0&-1&0 \\
        -2&4&3&-5&3&1&-2&0 \\ -2&3&3&-4&2&1&-2&1 \\ -2&2&3&-3&1&1&-1&0 \\
        -2&1&2&-2&1&1&-1&0 \\ -1&0&1&-1&1&0&0&0
        \end{array}\right)$,\\
  $\left(\begin{array}{cccccccc} -1&1&2&-2&0&1&0&0 \\ -2&2&2&-2&-1&2&0&0 \\ -3&2&4&-3&-1&2&0&0 \\
        -4&3&5&-4&-2&4&-1&1 \\ -3&2&4&-3&-2&3&0&1 \\
        -2&1&3&-2&-2&3&0&0 \\ -1&1&2&-2&-1&2&0&0 \\ -1&0&1&-1&0&1&0&0 \end{array}\right)$,
 $\left(\begin{array}{cccccccc} 0&-2&0&1&0&0&0&-1 \\ 1&-2&-1&1&1&-1&0&-1 \\
        1&-3&0&1&1&-1&0&-2 \\ 2&-4&-1&1&2&-1&0&-3 \\
        2&-3&-1&1&1&0&-1&-2 \\ 2&-3&-1&1&1&0&-1&-1 \\
        1&-2&-1&1&1&0&-1&-1 \\ 1&-1&-1&1&0&0&0&-1 \end{array}\right)$,
  $\left(\begin{array}{cccccccc} 0&-1&0&1&-1&1&0&-2 \\ 1&0&-1&1&-1&1&0&-3 \\ 1&-1&0&1&-1&1&0&-4 \\
        2&-1&-1&1&-1&2&0&-6 \\ 2&-1&-1&1&-1&2&-1&-4 \\
        2&-1&-1&1&-1&2&-1&-3 \\ 1&-1&-1&1&0&1&-1&-2 \\
        1&-1&-1&1&0&0&0&-1 \end{array}\right)$,\\
  $\left(\begin{array}{cccccccc} 0&1&0&-1&1&-1&0&2 \\ -1&0&1&-1&1&-1&0&3 \\ -1&1&0&-1&1&-1&0&4 \\
        -2&1&1&-1&1&-2&0&6 \\ -2&1&1&-1&1&-2&1&4 \\ -2&1&1&-1&1&-2&1&3 \\
        -1&1&1&-1&0&-1&1&2 \\ -1&1&1&-1&0&0&0&1 \end{array}\right)$,
  $\left(\begin{array}{cccccccc} 0&2&0&-1&0&0&0&1 \\ -1&2&1&-1&-1&1&0&1 \\ -1&3&0&-1&-1&1&0&2 \\
        -2&4&1&-1&-2&1&0&3 \\ -2&3&1&-1&-1&0&1&2 \\ -2&3&1&-1&-1&0&1&1 \\
        -1&2&1&-1&-1&0&1&1 \\ -1&1&1&-1&0&0&0&1 \end{array}\right)$,
  $\left(\begin{array}{cccccccc} 1&-1&-2&2&0&-1&0&0 \\ 2&-2&-2&2&1&-2&0&0 \\ 3&-2&-4&3&1&-2&0&0 \\
        4&-3&-5&4&2&-4&1&-1 \\ 3&-2&-4&3&2&-3&0&-1 \\
        2&-1&-3&2&2&-3&0&0 \\ 1&-1&-2&2&1&-2&0&0 \\ 1&0&-1&1&0&-1&0&0
        \end{array}\right)$,\\
  $\left(\begin{array}{cccccccc} 1&-1&-1&1&-1&0&1&0 \\ 1&-2&-1&2&-1&-1&1&0 \\
        2&-3&-2&3&-2&0&1&0 \\ 2&-4&-3&5&-3&-1&2&0 \\
        2&-3&-3&4&-2&-1&2&-1 \\ 2&-2&-3&3&-1&-1&1&0 \\
        2&-1&-2&2&-1&-1&1&0 \\ 1&0&-1&1&-1&0&0&0 \end{array}\right)$,
  $\left(\begin{array}{cccccccc} 1&-1&0&0&1&-1&0&1 \\ 0&-1&0&0&2&-2&0&2 \\ 0&-2&1&0&2&-2&0&2 \\
        0&-3&0&1&3&-3&0&3 \\ 0&-2&0&0&3&-2&0&2 \\ 0&-2&0&0&2&-1&0&2 \\
        0&-1&0&0&1&-1&1&1 \\ 0&0&0&0&0&0&0&1 \end{array}\right)$,
  $\left(\begin{array}{cccccccc} 1&-1&0&1&0&-1&-1&1 \\ 2&-1&-1&1&1&-2&-1&2 \\ 2&-1&0&1&1&-3&-1&2 \\
        4&-2&-1&2&1&-4&-1&3 \\ 3&-2&-1&2&1&-3&-1&2 \\
        2&-2&0&1&1&-2&-1&2 \\ 1&-1&0&1&0&-1&-1&2 \\ 1&0&0&0&0&0&-1&1
        \end{array}\right)$,\\
 $\left(\begin{array}{cccccccc} 1&0&-2&2&-1&0&0&-1 \\ 2&0&-2&2&-1&0&0&-2 \\
        3&0&-4&3&-1&0&0&-2 \\ 4&0&-5&4&-1&-1&1&-4 \\
        3&0&-4&3&0&-1&0&-3 \\ 2&1&-3&2&0&-1&0&-2 \\ 1&0&-2&2&0&-1&0&-1 \\
        1&0&-1&1&0&-1&0&0 \end{array}\right)$,
  $\left(\begin{array}{cccccccc} 1&0&-1&0&1&0&-2&2 \\ 1&0&0&-1&2&-1&-2&3 \\ 2&0&-1&-1&2&0&-3&4 \\
        3&1&-1&-2&3&-1&-4&6 \\ 2&1&-1&-1&2&-1&-3&5 \\
        1&1&-1&-1&2&-1&-2&4 \\ 0&1&0&-1&1&0&-2&3 \\ 0&1&0&-1&1&0&-1&1 \end{array}\right)$,
  $\left(\begin{array}{cccccccc} 1&0&-1&0&2&-2&0&0 \\ 1&-1&-1&0&3&-2&-1&1 \\
        1&0&-2&0&4&-3&-1&1 \\ 1&-1&-2&0&6&-4&-2&1 \\
        1&0&-1&-1&5&-3&-2&1 \\ 0&0&0&-1&4&-3&-1&1 \\ 0&0&0&-1&3&-2&0&0 \\
        0&0&0&0&1&-1&0&0 \end{array}\right)$,\\
  $\left(\begin{array}{cccccccc} 1&0&-1&1&-2&1&1&-1 \\ 1&0&-1&2&-3&1&1&-2 \\ 2&-1&-2&3&-4&2&1&-2 \\
        2&-1&-3&5&-6&2&2&-3 \\ 2&-1&-3&4&-4&1&2&-3 \\
        2&0&-3&3&-3&1&1&-2 \\ 2&0&-2&2&-2&0&1&-1 \\ 1&0&-1&1&-1&0&0&0 \end{array}\right)$,
  $\left(\begin{array}{cccccccc} 1&0&0&0&0&-1&0&0 \\ 2&-1&-1&1&0&-1&-1&1 \\ 2&0&-1&1&0&-2&-1&1 \\
        3&-1&-1&2&-1&-2&-2&2 \\ 3&-1&-1&2&-1&-2&-1&1 \\
        2&-1&0&1&0&-2&-1&1 \\ 2&0&0&0&0&-1&-1&1 \\ 1&0&0&0&0&0&-1&0 \end{array}\right)$,
  $\left(\begin{array}{cccccccc} -2&0&1&0&0&-1&1&1 \\ -2&-1&1&0&0&-1&2&1 \\ -3&0&1&0&0&-2&3&1 \\
        -4&-1&2&0&0&-3&5&1 \\ -4&-1&2&0&0&-2&4&0 \\ -3&-1&2&0&0&-2&3&0 \\
        -2&-1&1&1&-1&-1&2&0 \\ -1&0&1&0&-1&0&1&0
        \end{array}\right)$,\\
  $\left(\begin{array}{cccccccc} 1&0&0&1&-1&0&-1&0 \\ 2&1&-1&1&-1&0&-1&0 \\ 2&1&0&1&-1&-1&-1&0 \\
        4&1&-1&2&-2&-1&-1&0 \\ 3&0&-1&2&-1&-1&-1&0 \\
        2&0&0&1&-1&0&-1&0 \\ 1&0&0&1&-1&0&-1&1 \\ 1&0&0&0&0&0&-1&1 \end{array}\right)$,
  $\left(\begin{array}{cccccccc} 1&1&-1&0&0&1&-2&1 \\ 1&2&0&-1&0&1&-2&1 \\ 2&2&-1&-1&0&2&-3&2 \\
        3&4&-1&-2&0&2&-4&3 \\ 2&3&-1&-1&0&1&-3&3 \\ 1&3&-1&-1&0&1&-2&2 \\
        0&2&0&-1&0&1&-2&2 \\ 0&1&0&-1&1&0&-1&1 \end{array}%
\right)$.

The character table of $G^{s_{67}}$:\\
${}_{
 \!\!\!_{
\!\!\!\!\!_{ \!\!\!\!\!_{  \small
 \noindent
%

\noindent \noindent where A = E(3)
  = (-1+ER(-3))/2 = b3,
B = E(9)$^7$,C = -E(9)$^4$-E(9)$^7$,D = E(9)$^4$.

The generators of $G^{s_{68}}$ are:\\
$\left(%
\right)$.

The character table of $G^{s_{68}}$:\\
\noindent %

\noindent \noindent where A = 3*E(3)$^2$
  = (-3-3*ER(-3))/2 = -3-3b3,
B = 6*E(3)$^2$
  = -3-3*ER(-3) = -3-3i3,
C = 9*E(3)$^2$
  = (-9-9*ER(-3))/2 = -9-9b3,
D = -E(3)$^2$
  = (1+ER(-3))/2 = 1+b3,
E = 2*E(3)$^2$
  = -1-ER(-3) = -1-i3,
F = E(3)+2*E(3)$^2$
  = (-3-ER(-3))/2 = -2-b3,
G = E(3)-E(3)$^2$
  = ER(-3) = i3.

The generators of $G^{s_{69}}$ are:\\
$\left(%
\right)$.

The character table of $G^{s_{69}}$:\\
${}_{
 \!\!\!\!\!_{
\!\!\!\!\!_{ \!\!\!\!\! _{  \small
 \noindent
%

\noindent \noindent where A = E(3)
  = (-1+ER(-3))/2 = b3,
B = E(9)$^7$,C = -E(9)$^4$-E(9)$^7$,D = E(9)$^4$.

The generators of $G^{s_{70}}$ are:\\
 $\left(%
\right)$.

The character table of $G^{s_{70}}$:\\
\noindent %

\noindent \noindent where A = 3*E(3)$^2$
  = (-3-3*ER(-3))/2 = -3-3b3,
B = 6*E(3)$^2$
  = -3-3*ER(-3) = -3-3i3,
C = 9*E(3)$^2$
  = (-9-9*ER(-3))/2 = -9-9b3,
D = -E(3)$^2$
  = (1+ER(-3))/2 = 1+b3,
E = 2*E(3)$^2$
  = -1-ER(-3) = -1-i3,
F = 2*E(3)+E(3)$^2$
  = (-3+ER(-3))/2 = -1+b3,
G = -E(3)+E(3)$^2$
  = -ER(-3) = -i3.

The generators of $G^{s_{71}}$ are:\\
$\left(%
\right)$.

The character table of $G^{s_{71}}$:\\
${}_{
 \!\!\!\!\!\!\!_{
\!\!\!\!\!\!\!_{ \!\!\!\!\!\!\!  _{  \small
 \noindent
%

\noindent \noindent where A = E(3)$^2$
  = (-1-ER(-3))/2 = -1-b3,
B = -E(9)$^7$,C = E(9)$^4$+E(9)$^7$,D = -E(9)$^4$.

The generators of $G^{s_{72}}$ are:\\
$\left(%
\right)$.

The character table of $G^{s_{72}}$:\\
${}_{
 \!\!\!\!\!\!\!_{
\!\!\!\!\!\!\!_{ \!\!\!\!\!\!\!  _{  \small
 \noindent
%

\noindent \noindent where A = E(3)$^2$
  = (-1-ER(-3))/2 = -1-b3,
B = -E(9)$^2$,C = E(9)$^2$+E(9)$^5$,D = -E(9)$^5$.

The generators of $G^{s_{73}}$ are:\\
$\left(%
\right)$.

The character table of $G^{s_{73}}$:\\
${}_{
 \!\!\!\!\!_{
\!\!\!\!\!_{ \!\!\!\!\!  _{  \small
 \noindent
%

\noindent \noindent where A = -E(3)
  = (1-ER(-3))/2 = -b3,
B = E(4)
  = ER(-1) = i,
C = E(12)$^7$.

The generators of $G^{s_{74}}$ are:\\
$\left(%
\right)$.

The character table of $G^{s_{74}}$:\\
${}_{
 \!\!\!\!\!_{
\!\!\!\!\!_{ \!\!\!\!\!  _{  \small
 \noindent
%

\noindent \noindent where A = -2*E(4)
  = -2*ER(-1) = -2i,
B = -4*E(4)
  = -4*ER(-1) = -4i,
C = -6*E(4)
  = -6*ER(-1) = -6i,
D = E(4)
  = ER(-1) = i,
E = -1+E(4)
  = -1+ER(-1) = -1+i.

The generators of $G^{s_{75}}$ are:\\
$\left(%
\right)$.

The character table of $G^{s_{75}}$:\\
\noindent %

\noindent \noindent where A = E(3)$^2$
  = (-1-ER(-3))/2 = -1-b3,
B = 2*E(3)$^2$
  = -1-ER(-3) = -1-i3,
C = 3*E(3)$^2$
  = (-3-3*ER(-3))/2 = -3-3b3,
D = 4*E(3)$^2$
  = -2-2*ER(-3) = -2-2i3,
E = 2*E(4)
  = 2*ER(-1) = 2i,
F = 2*E(12)$^{11}$,G = 4*E(4)
  = 4*ER(-1) = 4i,
H = 4*E(12)$^{11}$,I = -E(4)
  = -ER(-1) = -i,
J = -E(12)$^{11}$,K = -1+E(4)
  = -1+ER(-1) = -1+i,
L = -E(12)$^8$+E(12)$^{11}$,M = -E(12)$^8$-E(12)$^{11}$.

The generators of $G^{s_{76}}$ are:\\
$\left(%
\right)$.

The character table of $G^{s_{76}}$:\\
\noindent %

\noindent \noindent where A = E(3)$^2$
  = (-1-ER(-3))/2 = -1-b3,
B = E(4)
  = ER(-1) = i,
C = E(12)$^{11}$,D = -2*E(4)
  = -2*ER(-1) = -2i,
E = -2*E(3)
  = 1-ER(-3) = 1-i3,
F = -2*E(12)$^7$.

The generators of $G^{s_{77}}$ are:\\
$\left(%
\right)$.

The character table of $G^{s_{77}}$:\\
${}_{
 \!\!\!\!\!_{
\!\!\!\!\!_{ \!\!\!\!\! _{  \small
 \noindent
%

\noindent \noindent where A = -E(3)$^2$
  = (1+ER(-3))/2 = 1+b3,
B = E(4)
  = ER(-1) = i,
C = E(12)$^{11}$,D = -2*E(4)
  = -2*ER(-1) = -2i,
E = 3*E(4)
  = 3*ER(-1) = 3i,
F = 2*E(3)
  = -1+ER(-3) = 2b3,
G = 3*E(3)
  = (-3+3*ER(-3))/2 = 3b3,
H = -2*E(12)$^{11}$,I = -3*E(12)$^{11}$.

The generators of $G^{s_{78}}$ are:\\
$\left(%
\right)$.

The character table of $G^{s_{78}}$:\\
\noindent %

\noindent \noindent where A = E(3)$^2$
  = (-1-ER(-3))/2 = -1-b3,
B = 2*E(3)$^2$
  = -1-ER(-3) = -1-i3,
C = 3*E(3)$^2$
  = (-3-3*ER(-3))/2 = -3-3b3,
D = 4*E(3)$^2$
  = -2-2*ER(-3) = -2-2i3,
E = 6*E(3)$^2$
  = -3-3*ER(-3) = -3-3i3,
F = 9*E(3)$^2$
  = (-9-9*ER(-3))/2 = -9-9b3.

The generators of $G^{s_{79}}$ are:\\
 $\left(%
\right)$.

The character table of $G^{s_{79}}$:\\
\noindent %

\noindent \noindent where A = E(4)
  = ER(-1) = i,
B = -E(3)$^2$
  = (1+ER(-3))/2 = 1+b3,
C = -E(12)$^7$,D = 1-E(4)
  = 1-ER(-1) = 1-i,
E = E(12)$^4$-E(12)$^7$,F = E(12)$^4$+E(12)$^7$,G = 2*E(3)$^2$
  = -1-ER(-3) = -1-i3,
H = -2*E(4)
  = -2*ER(-1) = -2i,
I = -2*E(12)$^7$,J = 3*E(3)
  = (-3+3*ER(-3))/2 = 3b3,
K = -4*E(4)
  = -4*ER(-1) = -4i,
L = -4*E(12)$^7$,M = 4*E(3)
  = -2+2*ER(-3) = 4b3.

The generators of $G^{s_{80}}$ are:\\
$\left(%
\right)$.

The character table of $G^{s_{80}}$:\\
\noindent %

\noindent \noindent where A = -4*E(4)
  = -4*ER(-1) = -4i,
B = -20*E(4)
  = -20*ER(-1) = -20i,
C = -36*E(4)
  = -36*ER(-1) = -36i,
D = -40*E(4)
  = -40*ER(-1) = -40i,
E = -64*E(4)
  = -64*ER(-1) = -64i,
F = 2*E(4)
  = 2*ER(-1) = 2i,
G = 6*E(4)
  = 6*ER(-1) = 6i,
H = -E(4)
  = -ER(-1) = -i,
I = 1+E(4)
  = 1+ER(-1) = 1+i,
J = 2-2*E(4)
  = 2-2*ER(-1) = 2-2i,
K = 6-6*E(4)
  = 6-6*ER(-1) = 6-6i,
L = 8*E(4)
  = 8*ER(-1) = 8i,
M = -4+4*E(4)
  = -4+4*ER(-1) = -4+4i.

The generators of $G^{s_{81}}$ are:\\
$\left(%
\right)$

The representatives of conjugacy classes of   $G^{s_{81}}$ are:\\
$\left(%
\right)$.

The character table of $G^{s_{81}}$:\\
\noindent %

\noindent \noindent where A = E(5)$^2$,B = E(5)$^4$,C = -E(4)
  = -ER(-1) = -i,
D = -E(20)$^{13}$,E = -E(20).

The generators of $G^{s_{82}}$ are:\\
$\left(%
\right)$.

The character table of $G^{s_{82}}$:\\
\noindent %

\noindent \noindent where A = E(3)$^2$
  = (-1-ER(-3))/2 = -1-b3,
B = 2*E(3)$^2$
  = -1-ER(-3) = -1-i3,
C = 3*E(3)$^2$
  = (-3-3*ER(-3))/2 = -3-3b3,
D = -E(4)
  = -ER(-1) = -i,
E = -E(12)$^{11}$,F = -2*E(4)
  = -2*ER(-1) = -2i,
G = -2*E(12)$^{11}$,H = -3*E(4)
  = -3*ER(-1) = -3i,
I = -3*E(12)$^{11}$.

The generators of $G^{s_{83}}$ are:\\
$\left(%
\right)$.

The character table of $G^{s_{83}}$:\\
\noindent %

\noindent \noindent where A = -E(4)
  = -ER(-1) = -i,
B = -2*E(4)
  = -2*ER(-1) = -2i,
C = -3*E(4)
  = -3*ER(-1) = -3i,
D = -4*E(4)
  = -4*ER(-1) = -4i,
E = -6*E(4)
  = -6*ER(-1) = -6i,
F = -9*E(4)
  = -9*ER(-1) = -9i.

The generators of $G^{s_{84}}$ are:\\
$\left(%
\right)$.

The character table of $G^{s_{84}}$:\\
${}_{
 \!\!\!\!\!\!\!_{
\!\!\!\!\!\!\!_{ \!\!\!\!\!\!\!  _{  \small
 \noindent
%

\noindent \noindent where A = E(3)$^2$
  = (-1-ER(-3))/2 = -1-b3,
B = -E(4)
  = -ER(-1) = -i,
C = -E(12)$^{11}$,D = 2*E(4)
  = 2*ER(-1) = 2i,
E = 2*E(12)$^{11}$,F = 2*E(3)$^2$
  = -1-ER(-3) = -1-i3.

The generators of $G^{s_{85}}$ are:\\
$\left(%
\right)$.

The character table of $G^{s_{85}}$:\\
${}_{
 \!\!\!\!\!\!\!\!\!\!\!\!\!_{
\!\!\!\!\!\!\!_{ \!\!\!\!\!\!\!  _{  \small
 \noindent
%

\noindent \noindent where A = E(7)$^2$,B = E(7)$^4$,C = E(7)$^6$.

The generators of $G^{s_{86}}$ are:\\
$\left(%
\right)$.

The character table of $G^{s_{86}}$:\\
${}_{
 \!\!\!\!\!\!\!_{
\!\!\!\!\!\!\!_{ \!\!\!\!\!\!\! _{  \small
 \noindent
%

\noindent \noindent where A = E(7)$^2$,B = E(7)$^4$,C = E(7)$^6$.

The generators of $G^{s_{87}}$ are:\\
$\left(%
\right)$.

The character table of $G^{s_{87}}$:\\
${}_{
 \!\!\!\!\!\!\!\!\!\!\!\!\!_{
\!\!\!\!\!\!\!_{ \!\!\!\!\!\!\!  _{  \small
 \noindent
%

\noindent \noindent where A = E(7)$^2$,B = E(7)$^4$,C = E(7)$^6$.

The generators of $G^{s_{88}}$ are:\\
$\left(%
\right)$.

The character table of $G^{s_{88}}$:\\
${}_{
 \!\!\!\!\!\!\!\!\!\!_{
\!\!\!\!\!\!\!_{ \!\!\!\!\!\!\!  _{  \small
 \noindent
%

\noindent \noindent where A = E(7)$^2$,B = E(7)$^4$,C = E(7)$^6$.

The generators of $G^{s_{89}}$ are:\\
$\left(%
\right)$.

The character table of $G^{s_{89}}$:\\
\noindent %

\noindent \noindent where A = E(3)$^2$
  = (-1-ER(-3))/2 = -1-b3,
B = 2*E(3)
  = -1+ER(-3) = 2b3,
C = 4*E(3)$^2$
  = -2-2*ER(-3) = -2-2i3.

The generators of $G^{s_{90}}$ are:\\
$\left(%
\right)$.

The character table of $G^{s_{90}}$:\\
${}_{
 \!\!\!\!\!\!\!_{
\!\!\!\!\!\!\!_{ \!\!\!\!\!\!\!  _{  \small
 \noindent
%

\noindent \noindent where A = -E(3)$^2$
  = (1+ER(-3))/2 = 1+b3,
B = -2*E(3)
  = 1-ER(-3) = 1-i3,
C = 3*E(3)
  = (-3+3*ER(-3))/2 = 3b3.

The generators of $G^{s_{91}}$ are:\\
$\left(%
\right)$.

The character table of $G^{s_{91}}$:\\
\noindent %

\noindent \noindent where A = E(4)
  = ER(-1) = i,
B = -1-E(4)
  = -1-ER(-1) = -1-i,
C = E(8)+E(8)$^3$
  = ER(-2) = i2,
D = E(8),E = -2*E(4)
  = -2*ER(-1) = -2i,
F = -4*E(4)
  = -4*ER(-1) = -4i,
G = E(8)-E(8)$^3$
  = ER(2) = r2,
H = 2*E(8),I = 3*E(4)
  = 3*ER(-1) = 3i,
J = 4*E(8).

The generators of $G^{s_{92}}$ are:\\
$\left(%
\right)$.

The character table of $G^{s_{92}}$:\\
${}_{
 \!\!\!\!\!\!_{
\!\!\!\!\!\!_{ \!\!\!\!\!\!  _{  \small
 \noindent
%
 }}}}$

\noindent \noindent where A = -E(4)
  = -ER(-1) = -i,
B = E(3)$^2$
  = (-1-ER(-3))/2 = -1-b3,
C = -E(12)$^{11}$,D = -E(8),E = E(24)$^{11}$,F = E(24)$^{19}$.

The generators of $G^{s_{93}}$ are:\\
$\left(%
\right)$.

The character table of $G^{s_{93}}$:\\
\noindent %

\noindent \noindent where A = -E(4)
  = -ER(-1) = -i,
B = -E(3)
  = (1-ER(-3))/2 = -b3,
C = -E(12)$^7$,D = 2*E(3)
  = -1+ER(-3) = 2b3,
E = -2*E(4)
  = -2*ER(-1) = -2i,
F = -2*E(12)$^7$.

The generators of $G^{s_{94}}$ are:\\
$\left(%
\right)$.

The character table of $G^{s_{94}}$:\\
\noindent %

\noindent \noindent where A = E(3)$^2$
  = (-1-ER(-3))/2 = -1-b3,
B = 2*E(3)$^2$
  = -1-ER(-3) = -1-i3,
C = 3*E(3)$^2$
  = (-3-3*ER(-3))/2 = -3-3b3,
D = 4*E(3)$^2$
  = -2-2*ER(-3) = -2-2i3,
E = 6*E(3)$^2$
  = -3-3*ER(-3) = -3-3i3.
\vskip 0.3cm

\noindent{\large\bf Acknowledgement}:\\ We would like to thank Prof.
N. Andruskiewitsch and Dr. F. Fantino for suggestions and help.

}
\end {document}